\documentclass[a4paper,10pt]{article}
\usepackage{amsfonts,amssymb,amsmath,amsthm}
\usepackage{mathrsfs}
\usepackage{cancel}
\usepackage{latexsym}

\usepackage[hmargin=2cm,vmargin=2.5cm]{geometry}

\title{The size of the nodal sets for the eigenfunctions of the smooth
laplacian}
\author{Demetrios A. Pliakis}

\begin{document}

\maketitle

\newcommand{\mangm}{\ensuremath{ ({\bf M}^n,g) }}
\newcommand{\manif}{\ensuremath{{\bf M}^n}}
\newcommand{\apgt} {\ {\raise-.5ex\hbox{$\buildrel>\over\sim$}}\ }
\newcommand{\aplt} {\ {\raise-.5ex\hbox{$\buildrel<\over\sim$}}\ }
\newcommand{\fo}{\ensuremath{\mbox{supp}}}
\newcommand{\sfe}[2]{\ensuremath{{\bf S}_{#1}^{#2}}}
\newcommand{\xw}[1]{\ensuremath{{\bf R}^{#1}}}
\newcommand{\hmi}{\ensuremath{\overline{{\bf R}_+}}}
\newcommand{\dian}[1]{\ensuremath{\underline{#1}}}
\newcommand{\kseks}{\begin{eqnarray}}
\newcommand{\tleks}{\end{eqnarray}}
\newcommand{\kseksw}{\begin{eqnarray*}}
\newcommand{\tleksw}{\end{eqnarray*}}
\newcommand{\lap}[1]{\ensuremath{\Delta_{#1}}}
\newcommand{\ei}{\ensuremath{\lambda}}
\newcommand{\eps}{\ensuremath{\epsilon}}
\newcommand{\vareps}{\ensuremath{\varepsilon}}
\newcommand{\eaf}[1]{\ensuremath{\sup_{#1}}}
\newcommand{\mkf}[1]{\ensuremath{\inf_{#1}}}
\newcommand{\eafor}[1]{\ensuremath{\sup_{I_{#1}}}}
\newcommand{\mkfor}[1]{\ensuremath{\inf_{I_{#1}}}}
\newcommand{\kl}{\ensuremath{\nabla}} 
\newcommand{\iq}{\ensuremath{\mbox{tr}}} 
\newcommand{\es}{\ensuremath{\nabla^2}}
\newcommand{\Sy}[1]{\ensuremath{\mbox{S}_{#1}}}
\newcommand{\klgs}{\ensuremath{\cancel{\nabla}}} 
\newcommand{\lapgs}{\ensuremath{\cancel{\Delta}}}
\newcommand{\varw}{\ensuremath{\upsilon}}
\newcommand{\ier}{\ensuremath{I_{r,\varepsilon}}}
\newcommand{\spapok}{\ensuremath{\cancel{\mbox{div}}}}
\newcommand{\spcurl}{\ensuremath{\cancel{\mbox{curl}}}}
\newcommand{\slih}{\ensuremath{\overline{h}}}
\newcommand{\slik}{\ensuremath{\overline{k}}}
\newcommand{\sli}[1]{\ensuremath{\overline{#1}}}
\newcommand{\sobp}{\ensuremath{\frac{n}{n-p}}}
\newcommand{\sosp}{\ensuremath{\frac{n-1}{n-p-1}}}
\newcommand{\isobp}{\ensuremath{ 1-\frac{p}{n} }}
\newcommand{\isosp}{\ensuremath{ 1-\frac{p}{n-1} }}
\newcommand{\dirfin}[2]{\ensuremath{{\mathbb D}^{1}[{#1};{\bf F};{#2}]}}
\newcommand{\Ha}[1]{\ensuremath{\mathscr{H}^{#1}}}

\newcommand{\harn}[1]{\ensuremath{{\mathrm H}_{#1}^2}}
\newcommand{\harz}[2]{\ensuremath{|#1|^{#2}}}
\newcommand{\hart}[1]{\ensuremath{\left|\frac{\nabla{#1}}{{#1}}\right|^2}}
\newcommand{\harl}[1]{\ensuremath{\left|\frac{\lap{}{#1}}{{#1}}\right|}}
\newcommand{\mpalla}[2]{\ensuremath{{\bf B}^{#1}_{#2}}}
\newcommand{\sphere}[2]{\ensuremath{{\bf S}^{#1}_{#2}}}
\newcommand{\mpallat}[1]{\ensuremath{\hat{B}_{#1}}}
\newcommand{\poir}[2]{\ensuremath{\dian{C}^{#1}_{#2},r_{#2}}}
\newcommand{\gener}[1]{\ensuremath{{\bf C}_{#1}}}
\newcommand{\gedi}[1]{\ensuremath{\mbox{d}_{#1}}} 
\newcommand{\arep}[1]{\ensuremath{{\bf C}^{#1} }}
\newcommand{\gmc}[1]{\ensuremath{{\bf D}^{#1} } }
\newcommand{\epik}[1]{\ensuremath{ {\bf P}^{#1} }}
\newcommand{\edra}[1]{\ensuremath{{\bf F}^{#1} }}
\newcommand{\brick}[1]{\ensuremath{{\bf P}^{#1} }}
\newcommand{\bricks}{\ensuremath{ {\bf P} }}
\newcommand{\face}[1]{\ensuremath{ {\bf F}^{#1} }}
\newcommand{\faces}{\ensuremath{ {\bf F} }}
\newcommand{\shell}[2]{\ensuremath{ {\bf K}^{#1}({#2}) }}
\newcommand{\shells}{\ensuremath{ {\bf K} }}
\newcommand{\bjes}[2]{\ensuremath{ ({#1};i_1,\dots,i_{#2}) }}
\newcommand{\fjes}[3]{\ensuremath{({#1};{i_1\dots i_{#2}};{#3})}}
\newcommand{\sjes}[3]{\ensuremath{({#1};{#2};{#3})}}
\newcommand{\smeg}[2]{\ensuremath{r;\eta_1,\dots,\eta_{#2};\epsilon_{#2}}}
\newcommand{\ire}{\ensuremath{I_{r,\varepsilon}}}
\newcommand{\gamp}[1]{\ensuremath{\gamma_{#1}(\phi)}}
\newcommand{\gamps}[1]{\ensuremath{\gamma_{#1}(\psi)}}
\newcommand{\gama}[1]{\ensuremath{\gamma_{#1}(\alpha)}}
\newcommand{\tube}{\ensuremath{{\bf T}}}
\newcommand{\sing}{\ensuremath{{\bf Sing}}}
\newcommand{\reg}{\ensuremath{{\bf Reg}}}
\newcommand{\qtube}{\ensuremath{{\bf Q}}}
\newcommand{\koh}[1]{\ensuremath{c_{#1}(H)}}
\newcommand{\hco}{\ensuremath{{\mathit D}}}
\newcommand{\nod}[2]{\ensuremath{{\bf N}_{#1}({#2})}}
\newcommand{\afa}{\ensuremath{A({\bf F})}}
\newcommand{\means}[1]{\ensuremath{\widehat{#1}}}
\newcommand{\ges}{\ensuremath{\cancel{\gamma}}}
\newcommand{\hs}{\ensuremath{\cancel{h}}}
\newcommand{\ks}{\ensuremath{\cancel{k}}}
\newcommand{\kas}{\ensuremath{\cancel{\kappa}}}
\newcommand{\lun}{\ensuremath{\ell,\nu}}
\newcommand{\ehr}{\ensuremath{\eps,\eta,\rho}}
\newcommand{\perio}[1]{\ensuremath{\overline{#1}}}
\newcommand{\errm}{\ensuremath{(m;\epsilon,\rho,r)}}
\newcommand{\echu}{\ensuremath{{\bf E}_{u,\chi}}}
\newcommand{\ogkos}[1]{\ensuremath{\mbox{vol}({#1})}} 
\newcommand{\diam}[1]{\ensuremath{\mbox{diam}({#1})}} 
\newcommand{\momb}[1]{\ensuremath{{#1}:{\bf P}\rightarrow {\bf R}}}
\newcommand{\mombp}[1]{\ensuremath{{#1}:\partial{\bf P}\rightarrow {\bf R}}}
\newcommand{\moms}[1]{\ensuremath{{#1}:{\bf K}\rightarrow {\bf R}}}
\newcommand{\momsp}[1]{\ensuremath{{#1}:\partial{\bf K}\rightarrow{\bf R}}}
\newcommand{\mwmf}[3]{\ensuremath{{#1}:{\bf F}^{#2}({#3})\rightarrow{\bf R}}}
\newcommand{\momwb}[1]{\ensuremath{{#1}:\partial{\bf W}\rightarrow {\bf R}}}
\newcommand{\momw}[1]{\ensuremath{{#1}:{\bf W}\rightarrow {\bf R}}}
\newcommand{\bolopo}[1]{\ensuremath{I_{{\bf B},#1}}}
\newcommand{\rs}{\ensuremath{\cancel{{\mathrm R}}}}
\newcommand{\duna}{\ensuremath{\Lambda}}
\newcommand{\cour}[1]{\ensuremath{{\bf C}_{#1}}}
\newcommand{\leda}{\ensuremath{{\bf L}}}
\newcommand{\kamps}{\ensuremath{{\bf K}}}
\newcommand{\kampsp}{\ensuremath{{\bf M}}}
\newcommand{\harm}[1]{\ensuremath{\widehat{#1}}}
\newcommand{\qwrio}{ \ensuremath{ {\bf W} }}
\newcommand{\dact}{\ensuremath{ {\bf A} }}
\newcommand{\ypf}[1]{\ensuremath{ \Upsilon_{#1} }}
\newcommand{\Dong}{\ensuremath{{\bf IU}}}
\newcommand{\dong}{\ensuremath{{\bf U}}}
\newcommand{\parad}[2]{\ensuremath{\frac{d#1}{d#2}}}
\newcommand{\paraf}[1]{\ensuremath{\frac{d{#1}}{dr}}}
\newcommand{\parat}[1]{\ensuremath{\frac{d^2{#1}}{dr^2}}}
\newcommand{\ric}{\ensuremath{{\mathrm Ric}}}
\newcommand{\R}{\ensuremath{{\mathrm R}}}
\newcommand{\Rm}{\ensuremath{{\mathrm Rm}}}
\newcommand{\B}{\ensuremath{{\mathrm  B}}}
\newcommand{\rie}{\ensuremath{{\mathrm Rm}}}
\newcommand{\se}{\ensuremath{{\mathrm K}}}
\newcommand{\olwf}{\ensuremath{\int_{{\bf F}}}}
\newcommand{\olws}{\ensuremath{\int_{{\bf K}}}}
\newcommand{\olwr}{\ensuremath{\int_{{\bf Reg}}}}
\newcommand{\olwt}{\ensuremath{\int_{{\bf T}}}}
\newcommand{\olwj}{\ensuremath{\int_{{\bf B}_j}}}

\newcommand{\olwsh}{\ensuremath{\int_{{\bf K}}}}
\newcommand{\olwfr}{\ensuremath{\int_{{\bf I}} }}
\newcommand{\invare}{\ensuremath{\int_{\mathcal{I}_{r,\varepsilon}}}}
\newcommand{\olwb}{\ensuremath{\int_{{\bf P}} }}
\newcommand{\olwbj}{\ensuremath{\int_{{ P}_j} }}
\newcommand{\olww}[1]{\ensuremath{ \int_{{\bf W}_{#1}} }}
\newcommand{\olwwr}[1]{\ensuremath{ \int_{{\bf A}(#1)}}}
\newcommand{\olwwm}[2]{\ensuremath{ \int_{{\bf W}_{#1}^{#2}} }}
\newcommand{\olwbn}[1]{\ensuremath{\int_{{\bf S}_{#1}} }}
\newcommand{\olwq}[1]{ \ensuremath{ \int_{{\bf W}_{#1}}} }
\newcommand{\oltube}[2]{\ensuremath{\int_{{\bf T}_{#1}({#2})} }}
\newtheorem{berniq}{Estimates}
\newtheorem{noc}[berniq]{Nodal volume}
\newtheorem{estimates}[berniq]{Theorem}
\newtheorem{nodvol}[berniq]{Theorem}
\newtheorem{hc}[berniq]{Theorem}
\newtheorem{alf}[berniq]{Theorem}
\newtheorem{ride}[berniq]{Lemma}
\newtheorem{lelo_in}[berniq]{Lemma}
\newtheorem{lo_in}[berniq]{Theorem}
\newtheorem{basint}[berniq]{Lemma}
\newtheorem{mcres}[berniq]{Lemma}
\newtheorem{hmcres}[berniq]{Lemma}
\newtheorem{ewfest}[berniq]{Lemma}
\newtheorem{harnshell}[berniq]{Lemma}
\newtheorem{mcraes}[berniq]{Lemma}
\newtheorem{int_by_parts}[berniq]{Lemma}
\newtheorem{noco}[berniq]{Lemma}

\paragraph{Description.}
A classical problem in physics and  geometry is the qualitative description 
of the spectrum of the laplacian on functions on a domain with Dirichlet 
boundary  conditions.  The eigenfunctions  determine the standing 
 vibrational modes of a drum  shaped as the given domain. 
It is expected that the 
eigenfunctions behave as polynomials of degree related to the order of the
harmonic. This is manifested for instance by the size of their  
nodal sets, i.e. the  set of zeroes. In terms of the vibrational motion of a
drum it consists of points that remain stationary and can be interpreted
as the locus of destructive interference of the waves with the boundary. 
In the late 18th century Chladni performed experiments with planar drums of
that revealed  the shape and size of the nodal set vary with the order of the harmonic
Similar questions can be asked for the laplacian acting on functions on 
a an arbitrary riemannian compact, closed manifold 

\paragraph{Statement and known results} Let 
\(\mangm, n\geq 2\) be an \(n\)-dimensional compact closed manifold equipped with a 
smooth riemannian metric \(g\). The Laplace-Beltrami operator  acting on
functions  \(\lap{g}\) on \(\mbox{M}^n\) is written in local coordinates,
\(g^{ij}=(g^{-1})_{ij}, g=\mbox{det}(g_{ij})\):
\[
\lap{g}\phi=\frac{1}{\sqrt{g}} \sum_{i,j}\frac{\partial}{\partial x_i}
\left( g^{ij}\sqrt{g}\frac{\partial \phi}{\partial x_j}  \right)
\]
The spectrum  of the laplacian is discrete
\(\mbox{spec}(\Delta)=\{\lambda_j\}_{j\in {\bf N}}\subset (0,\infty)\)
consists of eigenvalues with 
eigenfunctions \(u_\ei,\, \ei\in\mbox{spec}(\Delta)\) are 
the solutions of the equation
\[
\lap{g}u_\ei+\ei u_\ei=0 
 \]
The nodal set of an eigenfunction $u_\ei$ is defined as 
\[
\nod{}{u_\ei}=\{ x \in M^n/ u_\ei(x)=0\} 
 \]
and it is proved that this set is a smooth submanifold outside a set of 
\((n-1)\)-Hausdorff measure zero.  Br\"unning [1979] and \cite{yse}  
 showed that for a smooth surface the nodal length has a lower bound 
as:  
\[
\Ha{1}(\nod{}{u_\ei}) \geq C\ei^{1/2} 
 \]
and then Yau conjectured (\cite{yse}, pr. no 74) 
that the Hausdorff measure of this set grows as 
\[
\Ha{n-1}(\nod{}{u_\ei})\sim \ei^{1/2} 
\]
which is evidently the case for the spherical harmonics on the sphere. 
\cite{df1} established this when \(\mangm\) is analytic. 
Furthermore they showed \cite{df2} that for smooth surfaces one has that 
\[
\Ha{1}(\nod{}{u_\ei})\leq C\ei^{3/4} 
 \]
which was obtained by \cite{do} by different methods. Note that 
\cite{na} proved that 
 \[
\Ha{1}(\nod{}{u_\ei})\leq C'\ei\log\ei 
\]
Hardt and Simon  in their excellent work \cite{hs} showed for $C^{1,1}$ 
metrics that 
\[
\Ha{n-1}(\nod{}{u_\ei})\leq C"\sqrt{\ei}^{\sqrt{\ei}} 
\]
\cite{jm} based on Donnely-Fefferman work  obtained bounds for the
size of  tubular neighbouhoods of nodal sets for real analytic metrics. 
Recently \cite{sz},  \cite{cm}, \cite{ma} 
came up with new lower estimates of the nodal volume. 
The first authors provide an interesting formula for the size of 
level sets comprising an integral over the regular level sets of the
eigenfunction.  

\paragraph{Short description of the method.} Our  method is based on the
construction of a cell decomposition of a riemannian manifold with the 
boundaries of cells consisting of pieces of geodesic spheres. 
This resembles the picture
suggested by the Huyghens principle  governing wave propagation. 
The cells are selected so  that we are allowed  to follow an 
inductive argument restricting the
eigenfunction on the faces and get precise growth estimates for it.
Recall that the area of geodesic spheres determines the absorption 
rate of a wave in the 
course of propagation in a medium. The faces are constructed through the
introduction of a cluster of points: these are thought   
as a collection of sources of spherical waves, they define the geodesic
spheres that provide the pieces forming the cells. 
The domains that are formed are called {\it geodesic pixels} 
while the faces are called {\it undulating fronts} 
and the geodesic sphere pieces called {\it elementary wave fronts}.
The pixel  size is arranged essentially by the mean curvature of
the elementary wave front.  The conjecture is then proved by an
inductive argument combined in the lower bound with the isoperimetric 
inequality and eigenvalue estimates. The upper bound is obtained by an
elaboration of the Dong formula.

\paragraph{Notation}
We introduce  some notation. First we introduce the localized energy of a 
geodesic pixel:
\[
{\mathbb E}^{(1)}(\R,\B;\zeta)=
\olwb \zeta^2 \left(|\kl \R|^2+|\B|^2\right)
\]
and 
\[
{\mathbb E}^{(0)}(\R;\zeta)=\olwb \zeta^2|\R|^3
\]
where we introduce the {\it Bach tensor}
\[
\B_{ijk}=\frac{1}{n-2}\mbox{curl}(\ric)_{ijk}-\frac{1}{2(n-1)(n-2)}\left(
g_{ij}\R_k-g_{ik}\R_j\right)
\]
The faces of the pixel are pieces of geodesic spheres called {\it
elementary wave fronts} and are denoted by 
\(\faces\) with mean curvature \(h\) and 
localized tension \(T^1(h;\vartheta)\) given for a smooth test function
\(\vartheta,\,\, \fo{\vartheta}\subset\faces\):
\[
\mathbb{T}^1(h;\vartheta)=\olwf \vartheta^2|\klgs h|^2,\qquad \mathbb{T}^0(h;\vartheta)=\olwf
\vartheta^2 h^2
\]
Then we introduce two numbers 
\[
r(\bricks)=\sup_{\zeta\in C_0^\infty(\bricks)}
\left(\frac{{\mathbb E}^{(1)}(\R,\B;\zeta)}{{\mathbb E}^0(\R;\zeta)}\right)
\]
and
\[
t(\faces)=\sup_{\theta\in C_0^\infty(\faces)}
\left(\frac{{\mathbb T}^1(h;\theta)}{{\mathbb T}^0(h;\theta)}\right)
\]
The first theorem comprises Harnack estimates on such sets  

\begin{estimates}

Let \( \mangm\) be a riemannian manifold with scalar curvature 
function \(\R\). Then in a pixel 
\(\bricks, \, r(\bricks)={\bjes{\ell}{k}}=\mu\) 
with generic face \(\faces, \, t(\faces)=\eta \) 
we have that for positive constants explicitly calculated
\begin{eqnarray*}
\mkf{\faces}|h|\leq c_{11}\inf_{\faces}|h|+\eps,\\
\mkf{\faces}|k|^2\leq c_{12}\inf_{\faces}|k|^2 +\eps
\end{eqnarray*}
where \(c_{11},c_{12}\) depend in an explicit way on \(\faces\).
Furthermore let  \(u\) be an eigenfunction  with 
eigenvalue \(\ei\) and \(\eps\) a regular value of \(u\)      
\[
 \bricks_\eps=\bricks\cap \{ |u_\ei|>\eps\},\qquad \faces_\eps=
\faces\cap \{ |u_\ei|>\eps\}
\] 
Then  for explicitly calculated positive  constants 
\( c_{2j}=c_{2j}(\bricks,\ei), c_{3j}=c_{3j}(\faces), j=0,1,2,3,4\)
we have that

\begin{eqnarray*}
\eaf{\bricks_\eps}|u|\leq c_{20} \eps\\ 
\eaf{\bricks_\eps}(|\nabla u|)\leq c_{21} \mkf{\bricks_\eps}|\kl u|+c_{24}\\
\eaf{\bricks_\eps}(|\kl^2 u|)\leq c_{22} \mkf{\bricks_\eps}|\kl^2 u|+ c_{25}\\
\eaf{\faces_\eps}|u|\leq c_{30} \eps\\
\eaf{\faces_\eps}(|\klgs u|)\leq c_{31} \mkf{\faces_\eps}|\klgs u|+c_{34}\\
\eaf{\faces_\eps}(|\klgs^2 u|)\leq c_{32} \mkf{\faces_\eps}|\klgs^2 u|+ c_{35}
\end{eqnarray*}

\end{estimates}

These estimates provide the basis for the inductive argument and the 
application of Dong's formula:
 
\begin{nodvol}
Let \(u\) be as above then 
\[
\Ha{n-1}(\nod{}{u_\ei})\sim\sqrt{\ei} 
\]
\end{nodvol}

\section{Decomposition of a manifold}

\subsection{Quantitative Huygens principle}

We introduce the  set of points, chosen to lie on the nodal set  
\[
\gener{0}=\{ \dian{C}_i^0\}_{i=1}^N \subset \manif
\]
and lying at distance  
\[
d(\dian{C}_i^0,\dian{C}_j^0)\sim \gedi{0}\sim\frac{1}{\sqrt{\ei}}
\]
This is possible since in every geodesic ball of radius
\(O(\ei^{-\frac12})\) there is always a  zero. 
Furthrmore we  define  geodesic balls of radii \(r_i\), 
\(\mpalla{n}{\poir{0}{i}}\) bounded by geodesic
spheres \( \sphere{n-1}{\poir{0}{i}} \). 
These balls are taken to overlap in \((n+1)\)-ples : we set as
\(\epik{0;i_1\dots i_k}\) the  overlap region  of the balls :
\[
\epik{0;i_1\dots i_k}=\bigcap_{j=1}^k \mpalla{n}{\poir{0}{i_j}}
\] 
for \( k=2,\dots,n+1 \). These domains are bounded by
spherical regions denoted as \(\edra{0;i_1\dots i_k;\ell}\), 
and have interior curvature data
\[
\ric^{0;i_1\dots i_k}
\]
and face second fundamental form
\[
\hs^{0;i_1\dots i_k;\ell},\,\,\,\, \ks^{0;i_1\dots i_k;\ell},
\]
We call these sets {\it geodesic pixels}.  

We produce generations of  such pixels after the introduction of new 
centers and arrive at the collection of pixels after generation \(j\): 
\[
\epik{j;i_1\dots i_k}
\]
Its boundary consists of the  elementary wave fronts 
\(\edra{j;i_1,\dots,i_k;\ell}\) and curvature data:
\[
\ric^{j;i_1\dots i_k},\,\,\, \hs^{j;i_1\dots i_k}, \,\,\,\,
\ks^{j;i_1\dots i_K},  
\]
It is written in the form for
\[
\partial\epik{j;i_1\dots i_k}=
\bigcup_{j',k',\ell'} \edra{j';i_1\dots i_{k'};\ell'}
\]
The faces  \( \edra{j';i_1\dots i_{k'};\ell'} \) are called
{\it elementary wave fronts} (EWF). Each pixel defines homothetic EWF
spanned by the tubular neighbouhoods of the elementary wave fronts:  
\[
\shell{\sjes{j}{k}{\ell}}{}=
\ier\times \edra{\ell;i_1\dots i_k},\quad \ier=
((1-\varepsilon)r,(1+\varepsilon)r)\
\]
We introduce the  localized tension  \({\mathbb T}(h;\vartheta)\) of an EWF
given for a  smooth test function 
\(\vartheta,\fo{\vartheta}\subset\edra{\ell;i_1\dots i_k}\):
\[
{\mathbb T}^j(h;\vartheta)=\olwf \vartheta|\klgs^j h|^2 
\]
where \(\hs\) is the mean curvature of the (EWF).  
We introduce two numbers
\[
r(\bricks)=\sup_{\zeta\in C_0^\infty(\bricks)}
\left(\frac{{\mathbb E}^{(1)}(\R,\B;\zeta)}{{\mathbb E}^0(\R;\zeta)}\right)
\] 
and
\[
t(\faces)=\sup_{\theta\in C_0^\infty(\faces)}   
\left(\frac{{\mathbb T}^1(h;\theta)}{{\mathbb T}^0(h;\theta)}\right)
\]
Let
\(\eta_\ell,\mu\) be positive constants.  We say that a pixel \(\bricks\) 
with boundary consisting of  EWF
\[
\partial\bricks=\bigcup_{\ell=1}^{n+1}\faces_\ell
\]
satisfies an \((\eta_\ell,\mu) \)  condition if 
\[
t(\faces_\ell)\leq \eta_\ell,\qquad r(\bricks)\leq \mu
\]
Let \((\eta_{j';i_1\dots i_{k'};\ell},\mu_j,\eps_{\ell})\) be
positive constants.  We assume that the geodesic pixels 
\(\brick{\bjes{j}{\ell}}\) are selected so that their EWF 
\(\face{\fjes{j'}{k'}{\ell'}} \)
satify \((\eta_{j';i_1\dots i_{k'}},\eps_{\ell})\) 
estimate while  \((\mu_{j},\eps_{\ell})\) are the
parameters in the  curvature estimate. This set of pixels is a subset of 
compact closure in the neighbourhood of the zero section in \(TM\). In the
sequel we will try to establish the equations defining this set, as
consequences of the above integral inequalities.

\paragraph{The structure of the  metric in  geodesic polar coordinates.}
Every geodesic pixel has center some \(\dian{C}^j_i\) obtained in the    
\(j\)-th generation of center selection. We consider geodesic coordinates
from the neighboring pixels. Therefore let  
\( \mpalla{n}{\poir{j}{i}}\) be a  geodesic ball centered at the point
\( \dian{C}^j_i\)  and introduce polar coordinates through Gauss lemma.  
The metric is written then as:
\[
 g=dr^2+\ges(r)
\]
where \( \ges(r)\) is a riemannian metric on the geodesic sphere
\( \sphere{n-1}{\poir{j}{i}} =\partial \mpalla{n}{\poir{j}{i}} \) 
with second fundamental  form and mean curvature respectively \( \ks,\hs \).
Accordingly  we have the first and second variation equations
for the metric \(\gamma\), if we denote the radial derivative by \(\paraf{}\)
 while  the angular ones by  \(;j\) while in this section the symbol
\(\klgs\) denotes collectively the angular derivatives:
\begin{subequations}\label{variations}
\begin{align}
\paraf{\gamma} &= 2k, \label{first}\\
\paraf{k_{ij}} &= -k_{im}k_j^m-\R_{0i0j}\label{second} 
\end{align}
\end{subequations}
Then setting:
\[
\cancel{\varw}=\frac12\log\left(\mbox{det}(\gamma)\right)
\]
\[ 
\sigma=R_{i0j0}k^{ij}+k_m^ik_j^mk_i^j, 
\]
\[
\kappa=\kappa_0 =|k|
\]
we infer  that 
\begin{equation}\label{morefvsv}
\begin{split}
\paraf{\cancel{\varw}}=h \\
 \kappa\paraf{\kappa}=\sigma  
\end{split}
\end{equation}
Based on Newton's identities for symmetric polynomials we get the 
following inequality for \(\sigma\):
 \[
(3-C_{n,3})h\kappa^3 -h^3 -|\rie|\kappa
\leq \sigma \leq \frac{C_{n,3}+3}{3}h\kappa^2-\frac12h^3+|\rie|\kappa
\]
or 
 \[
-(C_{n,3}+2n)\kappa^3 -|\rie|\kappa
\leq \sigma \leq \frac{C_{n,3}+3n}{3}\kappa^3+|\rie|\kappa
\]

\paragraph{The structure equations}

The Gauss equations that relate the curvature of 
\(\gamma, \perio{R}\) to the ambient curvature
\begin{subequations}\label{gauss}
\begin{align}
\perio{R}_{imjn}+(k_{ij}k_{mn}-k_{in}k_{jm}) &= \R_{imjn},\label{gauss1} \\
\perio{R}_{ij}+k_{ij}h-k_{im}k^m_j & = \R_{ij}, \label{gauss2}\\
\perio{R}+h^2-k^2 = \R-\R_{00},\label{gauss3}
\end{align}
\end{subequations}
The Codazzi equations  
\kseksw
\klgs_ik_{jm}-\klgs_jk_{im}= R_{m0ij}, 
\tleksw
constitute a Hodge system:
\begin{equation}\label{hodge}
\begin{split}
\spcurl(k)_{ijm}=R_{moij}\\
\spapok(k)_i-\kl_i h=R_{0i}
\end{split}
\end{equation}
where 
\[
\spcurl(U)_{ijk}=\klgs_kU_{ij}-\klgs_jU_{ik},\,\,\,\,\,\spapok(U)_i=\klgs_jU^j_i
\]

\subsection{The second fundamental form and the mean curvature of the fronts}

\subsubsection{Harnack on the slice}

The Gauss-Codazzi equation for the spherical front gives that
\begin{equation}
\begin{split}
|\klgs h|^2 & = \spapok(k)_ih_i- R_{0i}h_i,\\
 |\spapok(k)|^2 & =\spapok(k)_ih_i+R_{0i}h_i
\end{split}
\end{equation}
Elaborating the preceding identities with Young inequality for sutiable 
\(p\) and obtain \(h\)-growth inequality in the spherical front domain
\( \faces =\face{\fjes{j}{k}{\ell}} \) and cut-off \(\vartheta\):
\[
\olwf |\spapok(\vartheta k)|^2 \leq \frac{\eta\eps^p}{p}\olwf
\vartheta^2|h|^{2p}+
\frac{(p-1) 
\ogkos{ \mathscr{F}_{k,\ell}}^{\frac{p+1}{2p}}}{p\eps^{\frac{p}{p-1}}}
\olwf |\spapok(\vartheta k)|^2
\]
or choosing \(\eps=\left(\frac{2(p-1)}{p}\right)^{1-\frac1p}
\ogkos{\mathscr{F}_{k,\ell}}^{\frac12-\frac{1}{p^2}}\) and get that
\[
\olwf |\spapok(\vartheta k)|^2\leq \eta\left(\frac2p\right)^p(p-1)^{p-1}
\ogkos{\mathscr{F}_{k,\ell}}^{\frac{p^2-1}{2p}}
\olwf \vartheta^2 h^{2p}
\]
 We recall Sobolev inequality from  \cite{lgmt}
for the case of EWF, \(\faces\subset\qwrio, r=\sosp\):
\[
\left( \olwf |U|^r \right)^{\frac1r} \leq C \olwf |\klgs U|+|h||U|
\] 
Starting from
\[
\left(\olwf |U|^{\ell r} \right)^{\frac1r} \leq \ell 
C\olwf |U|^{\ell-1} |\klgs U|+|h||U|^\ell
\]
and applying H\"older's inequality

\begin{equation}\label{sobsli}
\begin{split}
\olwf |U|^{\ell-1}|\klgs U|& \leq 
\left(\ogkos{\faces}\right)^{1-\frac1p}\left(\olwf 
|U|^{p(\ell-1)}\right)^{\frac1p}
\left(\olwf |\klgs U|^{\frac{p}{p-1}}\right)^{1-\frac1p} \\
\olwf|h||U|^\ell & \leq \left(\olwf
|U|^{q\ell}\right)^{1/q}\left(\olwf|h|^{\frac{q}{q-1}}\right)^{1-\frac1q}
\end{split}
\end{equation} 
where  \(q<2\). The inequalities on the slice
\[
\olwf |\klgs U|^2\leq C\olwf |\spapok(U)|^2+|\spcurl(U)|^2+|\rie||U|^2
\]
read for the localization second fundamental form, \(U=\zeta |k|\):
\[
\olwf \zeta^2|\klgs k|^2\leq C\olwf |\zeta|^2|\ric|+\zeta|\klgs h|^2
+\left(|\klgs \zeta|^2+|\rie|\zeta^2\right)|k|^2
\]
The last term in the right hand side leads after the application of
Young's  inequality combined with Sobolev's inequality to
\[
\olwf \zeta^2|\klgs k|^2\leq C\olwf |\zeta|^2|\ric|+\zeta|\klgs h|^2
+\olwf\left(|\klgs \zeta|^2+|\rie|\right)^{\frac{n-1}{2}}
\]
We consider now the slice regions \(\faces_j\) determined 
by the tension energy through the sequence of constants \(\{\eta_j\}\),
in which \(t(\faces_j)=\eta_j\). The Harnack inequality is proved through Moser iteration 
on the domain with smooth boundary  \(\qwrio \subset \faces\)   
obtained by smoothing out the boundary of \(\faces_j\).
Therefore we exhaust the domain through the harmonic approximation of
the face defining function \(F,\harm{F}_0\):
\[
\qwrio_j(\eta)=\{\dian{x}\in{\bf F}:
(\theta-\theta^{j+1})\eta \leq
|\harm{F}_0(\dian{x})|\leq (1-\theta+\theta^j)\eta\}
\]
and then Harnack inequality takes the form:
\[
\eaf{{\bf W}_\infty}|h|\leq \hco(\eta,\eta_\ell) 
\mkf{{\bf W}_\infty}|h|
\]
where the quantity \(\hco(\eta,\eta_\ell)>0\) is calculated in the appendix.

\paragraph{The growth of the tension integral}

The preceding estimates necessitate the derivation of 
radial growth  estimates for the tension integral 
\[
T={\mathbb T}^1(h;\vartheta)
\]
We start differentiating and  proceed with the application of the   
structural equations. We obtain the differentail inequality 
\kseksw
\paraf{T}\leq 2\left(\olwf \vartheta^2|h||\klgs h|^2 +\olwf
|\klgs h|^2\vartheta\left|\paraf{\vartheta}\right|+
\olwf \vartheta^2|h^i|\left|\paraf{h_i}\right|\right)
\tleksw
We obtain that 
\[
\left|\paraf{h_{i}}\right|\leq 2|k||\klgs k|+|\klgs\R_{00}|
\]
Therefore we have that:
\[
\paraf{T}\leq \left(\sup_{\edra{}}|h|\right)T+
\left(\sup_{\edra{}}|k|\right)
\left(\olwf |\klgs h|^2\right)^{1/2}\left(\olwf |\klgs k|^2\right)^{1/2}
\]
The last term is majorised by 
\[
\left(\olwf |\klgs h|^2\right)^{1/2}\left(\olwf |\klgs k|^2\right)^{1/2}\leq
CT+\olwf \left(\zeta^2|\ric|+|\rie||\klgs \zeta|^2\right)
\]
Furthermore we have that:
\[
\sup_{\edra{}}|k|\leq c_1T^{\frac12}+c_2
\]
We  select cutoffs satisfying for some \(\eta>0\)
\[
\paraf{|\klgs\vartheta|}\leq \eta|\klgs \vartheta| 
\]
We arrive at the inequality 
\[
\paraf{T}\leq c_1T^{3/2}+c_2, \quad c_i=c_i(\edra{},\eps),\,
i=1,2
\]
We conclude through the use of Young's inequality:
\[
\paraf{T}\leq c\left(T^2+1\right)
\]
and hence for \(r\leq \frac{9\pi}{20C}\):
\[
\left|T(r(1+\eps))-T(r(1-\eps))\right|\leq 
 \tan(Cr)\leq Cr
\]

\subsubsection{Radial Harnack estimates} 

In this section we derive the radial variation of the curvature quantities 
in the radial interval \(\ier=((1-\varepsilon) r,(1+\varepsilon)r)\)
 relative to  a given value of this quatities  at \(r\).
We write \eqref{second} in  contracted form:
\kseks
\paraf{h}=-k^2-\R_{00} 
\tleks
Fisrt we have that 
\[
\paraf{h} \leq -\frac{1}{n-1} h^2-\R_{00}\Leftrightarrow h^2\leq
-(n-1)\paraf{h}-(n-1)\R_{00}
\] 
We derive the differential inequality for \(k\):
\[
\paraf{|k|}\leq |k|^2+|\rie|\leq (|k|^2+1)(|\rie|+1)
\]
that is written after integration and elementary trigonometry as:
\[
\frac{||k(r(1+\vareps))|-|k(r(1-\vareps))||}{1+|k|(r(1+\vareps))|k|(r(1-\vareps))}
\leq \tan(Cr)
\]
for 
\[
\max_{\ire}\left(C(\vareps)\olws |\rie|\right), 
\max_{\ire}\left(|k(r(1+\vareps))|\cdot|k(r(1-\vareps)|\right)\leq C
\]
For \(r\leq\frac{9\pi}{20C}\):
\[
||k(r(1+\vareps))|-|k(r(1-\vareps))||\leq C'r
\]
The estimates for \(|\kl k|,|\kl^2k|\) follow from the elementary
differential inequalities:
\[
\paraf{|\kl k|}\leq |k||\kl k|+(|\kl \rie|+|\rie||k|)
\]
\[
\paraf{|\kl^2 k|}\leq C\left(|k||\kl^2 k|+|\kl k|^2+|\rie||\kl k|+
|\kl \rie||k|+|\kl^2\rie|\right)
\]
This lead to the estimates:
\[
|\kl k|\leq \left[Cr^2\left(\olws|\rie|\right)+
\olws|\kl \rie|\right]e^{\frac{C^2r^2}{2}}
\]
Similarly
\[
|\kl^2 k|\leq \left[Cr^2\left(\olws|\rie|\right)+\olws|\kl \rie|\right]e^{\frac{C^2r^2}{2}}
\]
Let \(\chi\) be a  cutoff supported in \(\fo{(\chi)}\subset  \ier\) such that 
\[
 \varepsilon |\chi'| +\varepsilon^2|\chi''|\leq C
\]
as well as 
\[
\mu_\varepsilon(r)=\invare \chi h(r',\dian{\xi})dr'
\]

\begin{mcres}
The following estimates hold for \(r  \leq \frac{1}{\sqrt{C_0}} \) and 
are relative to a fixed value of \(h(r)\):
\kseksw
|\mu_\varepsilon| \leq (2\varepsilon r)^{1/2}
\left(|h|^{1/2} + (2\varepsilon r|C_0^0|)^{1/2}\right)\\
|\klgs \mu_\varepsilon|\leq  2\vareps r
\left(|\klgs h|+ 2(1+\vareps)r c_1(|\rie|,|\kl \rie|)
\right)\\
|\klgs^2\mu_\varepsilon|\leq  (2\varepsilon r)
\left(|\klgs^2h|+2r(1+\vareps)c_2(|\rie|,|\kl \rie|,|\kl^2\rie|)\right)
\tleksw
\end{mcres}

\paragraph{The estimate of \(\mu_\varepsilon\)} We start applying  
and Cauchy-Schwarz in the end
\[
|\mu_\varepsilon|\leq (2\varepsilon r)^{1/2}
\left(\invare h^2\right)^{1/2}\leq ((n-1)\vareps r)^{1/2}
\left(|h(r(1-\vareps))-h(r(1+\vareps))|^{1/2}+(2\varepsilon C_0r)^{1/2}\right)
\leq (n-1)^{\frac12}r\vareps(|k|+C_0) 
\]

\paragraph{Estimate of \(|\klgs \mu_\varepsilon|, |\klgs^2\mu_\vareps|\)}

We commence with 
the integration by parts for \(U:\ire\rightarrow \xw{}\):
\[
\left|\invare U \right|=\left|\invare \paraf{r}U\right|\leq 
\invare \paraf{r}|U|\leq 
2\varepsilon r|U|- \invare r\paraf{|U|}
\]
and apply this to
\[
\paraf{|\klgs h|}\leq C_1\left(|k||\klgs k|+|\klgs\R_{00}|\right)
\]
\[ 
\paraf{|\klgs^2h|}\leq C_2\left(|k||\klgs^2k|+|\klgs k|^2+ 
|\klgs^2\R_{00}|\right)
\]
and obtain the desired estimates.

\paragraph{The radial-slice estimates.} 
We start recalling the differentiation identity:
\[
\paraf{}\left(\olwf U\right) =\olwf \paraf{U}+hU
\]
We introduce a cutoff function \(\harm{\chi}\) supported in the 
shell \( \fo{\harm{\chi}}\subset \shells\):
\kseksw
\shells= \ire \times \faces,\\
\ire=\left((1-\varepsilon)r,(1+\varepsilon)r\right)
\tleksw
Therefore applying Sobolev inequality on the slice, H\"older and standard 
elliptic estimates we derive the differential inequalities 
for the quantities
\[
\paraf{U_j}\leq a U_j+b_j, \qquad \paraf{H_j}\leq a\sqrt{U_j}H_j+c_j
\]
where
\begin{subequations}
\begin{align}
U_j&= \olwf |\klgs_jk|^2,\qquad H_j=\olwf |\klgs_jh|^2\\
a &= \hco(\eta)\left[\olwf |k|^2+h^2\right]^{1/2}, b_0 = \olwf |\rie|,\\
 b_1 &= \olwf |\kl\rie|^2+U_0\olwf|\rie|^2,\\
b_2 &= U_1\left[U_1+U_0\olwf|\rie|^2\right]+
\left(\olwf|\kl \rie|^2\right)U_0+\olwf|\kl^2\rie|^2
\end{align}
\end{subequations}
The diffrential equations  for \(U_j,H_j\) are of the form
\[
\paraf{U}\leq aU+b
\]
 we obtain for \(|\varrho|\leq \vareps\):
\[
U(r(1+\varrho))\leq e^{r\int_0^{\eps}a(1+r')dr'}\left(U(r)+r\int_0^{\varrho}
b(1+r')e^{-r'\int_0^{\varrho'}a(1+\varrho')d\varrho'}\right)
\]

\begin{hmcres}

The following hold true in \(\ire\) for \(i,j=1,2\) and constants
\(c_{ij}=c_{ij}(r;\eta_1,\dots,\eta_j;\eps_j)>0\)
\[
U_{j,\pm}(r(1\pm \eps))\leq c_{1j} U_{j,\pm}(r),\qquad
H_{j,\pm}(r(1\pm \eps))\leq c_{2j} H_{j,\pm}(r)
\]
\end{hmcres}

\subsection{Selection of the pixels through the front estimates}

\subsubsection{The basic ansatz}
We introduce the notation for a smooth domain \(\qwrio\) in a  riemannian manifold 
equipped with riemannian volume \(dv\):
\[
{\cal D}^j(U;\qwrio)=\int_\qwrio|\kl^jU|^2dv
\]
We will derive the equation satisfied by \(u_\ei\) near (EWF) i.e. in a 
spherical shell of the form 
\[
\shells= \ire \times \faces= ((1-\varepsilon)r,(1+\varepsilon)r)\times
\faces
\] 
for some \(\faces=\face{\fjes{j}{k}{\ell}}\)
with coordinates \( (r,\dian{\theta}) \) and volume 
\( dr d\sigma = \sqrt{\gamma}drd\theta\), while  \(h_\ell\) stands for the 
mean curvature of the shell. 

We  select the front so that its mean
curvature is controlled  by the eigenfunction growth near it.  
The eigenfunction equation  is written near (EWF)  in the form: 
\[
\parad{^2u_\ei}{r^2}+h\paraf{u_\ei}+\lap{\gamma}u_\ei=-\ei u_\ei
\]
We make the ansatz for the parameter \(\beta\) to be determined and the 
smooth functions \(\alpha,\phi\):  
\[
u_\ei(r,\theta)=
A\left(r,\dian{\theta}\right)
\sin\left(\beta_1\phi(r,\dian{\theta})\right),\,\,\,\,
\]
and  arrange the (EWF) so that near the front \(A\) is smooth and 
positive which implies that \(\phi\) inherits the unique continuation
property from \(u_\ei\).  Note also that we shoule inspect also level sets
of
the form \(\frac{k\pi}{\beta-1}\) for suitable \(k\)'s.
The equation then splits as:
\begin{equation}\label{ansatz}
\begin{split}
\left[\parat{A}+h\paraf{A}
-\beta_1^2\left(\left(\paraf{\phi}\right)^2+
|\klgs\phi|^2\right)A+
\left(\lapgs A+\ei A\right)\right]\sin(\beta_1\phi)+\\
+\beta_1\left[\left(2\paraf{A}+hA\right) 
\paraf{\phi}
+A\parat{\phi}+\left(2\klgs A\cdot\klgs \phi+
A\lapgs\phi\right)\right]\cos(\beta_1\phi)=0
\end{split}
\end{equation}

\subsubsection{Determination of local phase and local amplitude}
We require that for a function \(\alpha\) and a parameter \(\beta_2\) 
to be  determined there holds in the interval \((r_0(1-r),r_0(1+r))\):
\[
\paraf{A}+\frac{h}{2}A=\frac{\beta_2}{2}\paraf{e^{\beta_2\alpha}}
e^{-\frac{\mu_\epsilon}{2}}
\]
with the conditions
\[
A(0,\theta)=1,\quad \alpha(0,\theta)=-\frac{2\log\beta_2}{\beta_2}
\]
We find 
\[
A(\sqrt{\ei}r,\theta)=e^\Lambda,\,\, 
\Lambda=-\frac12\mu_\varepsilon+\beta_2\alpha
\]
Therefore we obtain the pair of first order equations:
\[
\left(\paraf{\phi}\right)^2+|\klgs\phi|^2=\frac{\ei}{\beta_1^2}-
\frac{1}{\beta_1^2}s_1, \quad
s_1=\beta_2\left(\lapgs \alpha +\parat{\alpha}\right)+
\beta_2^2\left[|\klgs\alpha|^2+\left(\paraf{\alpha}\right)^2\right]-
\beta_2\klgs\mu\cdot\klgs\alpha+\kamps
\]
if 
\[
\kamps=\frac14\left(2|k|^2+2R_{00}-h^2\right)-\frac12\lapgs\mu+
\frac14|\klgs\mu|^2
\]
The other term gives  
\[
\beta_2
\left(\paraf{\phi}\paraf{\alpha}+2\klgs\phi\cdot\klgs\alpha\right)=-s_2, \quad 
s_2=\parat{\phi}+\lapgs\phi-\klgs\mu\cdot \klgs \phi
\]
We choose \(\beta_1=\sqrt{\ei}=\frac{1}{\beta_2}\). 
The solution of these equations is achieved through the method of 
characteristics.   We set 
\begin{subequations}\label{auxifu}
\begin{align}
{\bf A}&=|\klgs\alpha|^2+\left(\paraf{\alpha}\right)^2\label{auxifu1} \\
{\bf P}& =|\klgs\phi|^2+\left(\paraf{\phi}\right)^2 \label{auxifu2}\\
\duna &=\frac2\ei\alpha +\mu \label{auxifu3}
\end{align}
\end{subequations} 
Then we will consider this system in the form 
\begin{subequations}\label{systemap}
\begin{align}
{\bf P}&= 1+
\frac{1}{\ei\sqrt{\ei}}\left[
\lapgs \alpha +\parat{\alpha}+
\frac{1}{\ei}{\bf A}-
\frac{1}{\sqrt{\ei}}\klgs\mu\cdot\klgs\alpha+\kamps\right]
\label{systemap1}\\
\paraf{\phi}\paraf{\alpha}+2\klgs\phi\cdot\klgs\alpha &= 
-\ei\left(\parat{\phi}+\lapgs\phi-\klgs\mu\cdot \klgs \phi\right)
\label{systemap2}
\end{align}
\end{subequations}
We derive estimates for \(\alpha,\phi\) in the form 
\begin{subequations}\label{systemap}
\begin{align}
-\lapgs \alpha -\parat{\alpha}-
\frac{1}{\sqrt{\ei}}\klgs\mu\cdot\klgs\alpha
& =\ei\sqrt{\ei}\left[1-{\bf P}\right]+\frac{1}{\ei}{\bf A}
+\kamps,
\label{systemap1}\\
-\parat{\phi}-\lapgs\phi+\klgs\Lambda\cdot \klgs \phi &=\frac1\ei
\paraf{\phi}\paraf{\alpha}
\label{systemap2}
\end{align}
\end{subequations}

\subsubsection{Equations for the higher derivatives of the phase and
amplitude functions}

Furthermore differentiating the equations \eqref{systemap} we obtain the necessary
equations for the higher derivatives of the phase function, 
\[
\Psi_{0,m}=\frac{d^m\phi}{dr^m},\,\, a_{0,m}=\frac{d^m\alpha}{dr^m}\,\,
\duna_{0,m}=\frac{d^m\duna}{dr^m}, \,\,
{\bf A}_m=\frac{d^m{\bf A}}{dr^m},\,\, {\bf P}_{0,m}=\frac{d^m{\bf P}}{dr^m}, 
\qquad m=1,2,3 
\]
and the angular derivatives 
\[
\Psi_{\gamma,m}=|\klgs^m\phi|^2,\,\, a_{\gamma,m}=|\klgs^m\alpha|^2
\,\,\duna_{\gamma,m}=|\klgs^m \duna|^2, \,\,
{\bf A}_{\gamma,m}=|\klgs^m{\bf A}|^2,\,\, {\bf P}_{\gamma,m}=
|\klgs^m{\bf P}|^2,  \qquad m=1,2,3 
\]

\paragraph{Radial derivatives of the phase function}

Specifically we obtain, applying the commutation rules, the following
equations:
\begin{subequations}\label{systemap2}
\begin{align}
\lapgs \Psi_{0,1}+ \parat{\Psi_{0,1}}- \klgs\duna \cdot
\klgs\Psi_{0,1} &={\bf  S}_{0,1} \label{systemap21}\\
\lapgs \Psi_{0,2}+ \parat{\Psi_{0,2}}- \klgs\duna \cdot \klgs
\Psi_{0,2} &={\bf S}_{0,2} \label{systemap23}\\
\lapgs \Psi_{0,3}+\parat{\Psi_{0,3}}- \klgs\duna \cdot \klgs \Psi_{0,3} &=
{\bf S}_{0,3} \label{systemap24}
\end{align}
\end{subequations}
where the ''source terms'' are the following:
\begin{subequations}\label{sources}
\begin{align}
{\bf S}_{0,1} &=R^0_{\,0}\Psi_1+\R^l_{\,0}\klgs_l\phi+\klgs\duna^1\cdot\klgs\phi
\label{source1}\\
{\bf S}_{0,2} &=2\R^l_{\,j}\klgs_j\Psi_{1}+\paraf{\R^l_{\,0}}\klgs_l\Psi_1
+\klgs\duna^2\cdot\klgs\phi
\label{source2}\\
{\bf S}_{0,3} &=-\left(\R^l_{\,0}\klgs_l\psi_2+ 2\R^l_{\,0j0;j} 
\klgs_l\Psi_1+2\R^l_{\,0j0}\klgs_l\Psi_1\klgs^j\Psi_1
+\klgs\duna^3\cdot\klgs\phi\right) 
\label{source3}
\end{align}
\end{subequations}
\paragraph{Higher angular derivatives of the phase} 
We obtain, applying the commutation rules, the following
equations:
\begin{subequations}\label{systemap4}
\begin{align}
\lapgs \Psi_{\gamma,1}+\parat{\Psi_{\gamma,1}}- \klgs\duna \cdot 
\klgs \Psi_{\gamma,1}&=
{\bf S}_{2,1} \label{systemap41}\\
\lapgs \Psi_{\gamma,2}+ \parat{\Psi_{\gamma,2}}- \klgs\duna \cdot \klgs \Psi_{\gamma,2} &=
{\bf S}_{2,2} \label{systemap42}
\end{align}
\end{subequations}
where the ''source terms'' are the following:
\begin{subequations}\label{sources}
\begin{align}
{\bf S}_{2,1} &= -\left(2\R^s_{\,00j}\klgs_j\phi+\R^l_{\,j}\klgs_j\phi\right)
\klgs_l\phi-(\tau^i_j\klgs_j\phi\klgs_i\phi+
\tau^i\klgs_i\klgs_j\phi\klgs_j\phi) 
\label{ssource1}\\
{\bf S}_{2,2} &= -\Rm*\klgs^2\phi*\klgs^2\phi-
\klgs\Rm*\klgs\phi*\klgs^2\phi-\Rm*\klgs\Psi_{0,1}*\klgs\phi+
\lapgs\tau^i\klgs_i\phi \label{ssource2}
\end{align}
\end{subequations}

\paragraph{Radial derivatives of the amplitude}
Set first that:
\[
{\bf V}=\frac1\ei{\bf A}+\ei^{3/2}\left(1-{\bf P}\right)+\kamps
\]
Specifically we obtain, applying the commutation rules, the following
equations:
\begin{subequations}\label{systemap2}
\begin{align}
-\lapgs a_{0,1}- \parat{a_{0,1}}+\frac{1}{\sqrt{\ei}}\klgs \mu\cdot
\klgs a_{0,1} &={\bf  S}_{3,1} +\paraf{{\bf V}}\label{systemap21}\\
-\lapgs a_{0,2}- \parat{a_{0,2}}+\frac{1}{\sqrt{\ei}}\klgs\mu \cdot \klgs
a_{0,2} &=
{\bf S}_{3,2}+\parat{{\bf V}} \label{systemap23}\\
-\lapgs a_{0,3}-\parat{a_{0,3}}+\frac{1}{\sqrt{\ei}}\klgs\mu \cdot 
\klgs a_{0,3} &= {\bf S}_{3,3}+\frac{d^3{\bf V}}{dr^3} \label{systemap24}
\end{align}
\end{subequations}
where the ''source terms'' are the following:
\begin{subequations}\label{sources}
\begin{align}
{\bf S}_{3,1} &=
R^0_{\,0}a_{0,1}+\R^l_{\,0}\klgs_l\alpha+\klgs\mu\cdot\klgs\alpha
\label{source1}\\
{\bf S}_{3,2} &=2\R^l_{\,j}\klgs_ja_{0,1}+\paraf{\R^l_{\,0}}\klgs_la_{0,1}
+\klgs\mu^{''}\cdot\klgs\alpha
\label{source2}\\
{\bf S}_{3,3} &=-\left(\R^l_{\,0}\klgs_la_{0,2}+ 2\R^l_{\,0j0;j} 
\klgs_la_{0,1}+2\R^l_{\,0j0}\klgs_la_{0,1}\klgs^j\Psi_{0,1}
+\klgs\mu^{'''}\cdot\klgs\alpha\right) 
\label{source3}
\end{align}
\end{subequations}

\paragraph{Higher angular derivatives of the amplitude} 
We obtain, applying the commutation rules, the following
equations:
\begin{subequations}\label{systemap4}
\begin{align}
\lapgs a_{\gamma,1}+\parat{a_{\gamma,1}}- \klgs\mu \cdot \klgs a_{\gamma,1} &=
{\bf S}_{4,1} \label{systemap41}\\
\lapgs a_{\gamma,2}+ \parat{a_{\gamma,2}}- \klgs\mu \cdot \klgs a_{\gamma,2} &=
{\bf S}_{4,2} \label{systemap42}
\end{align}
\end{subequations}
where the ''source terms'' are the following:
\begin{subequations}\label{sources}
\begin{align}
{\bf S}_{4,1} &= -\left(2\R^s_{\,00j}\klgs_j\alpha+\R^l_{\,j}\klgs_j\alpha
\right)\klgs_l\alpha-(\tau^i_j\klgs_j\alpha\klgs_i\alpha+
\tau^i\klgs_i\klgs_j\alpha\klgs_j\alpha) 
\label{ssource1}\\
{\bf S}_{4,2} &= -\Rm*\klgs^2\alpha*\klgs^2\alpha-
\klgs\Rm*\klgs\alpha*\klgs^2\alpha-\Rm*\klgs a_{0,1}*\klgs\alpha+
\lapgs\tau^i\klgs_i\alpha \label{ssource2}
\end{align}
\end{subequations}

\subsubsection{Estimates }
We recall the variation identities:
\begin{subequations}\label{raddif}
\begin{align}
 \olwf v\paraf{v}& =\frac12\paraf{}\left(\olwf v^2 \right)-\frac12\olwf hv^2 
\label{raddif1}\\
  \olwf v\parat{v}& = \frac12\parat{}\left(\olwf v^2\right)
-\olwf\left(\paraf{v}\right)^2-\paraf{}\left(\olwf hv^2\right)
\label{raddif2} \\
\paraf{h} & =-|k|^2-R_{00}
\end{align}
\end{subequations}
\eqref{systemap1}  gives the pointwise estimate  provided \(\alpha\geq 0\):
\kseks\label{pointamp}
-\alpha\lapgs \alpha-\alpha\parat{\alpha}
-\frac{1}{\ei}{\bf A}
\alpha +\frac{1}{\sqrt{\ei}}(\klgs \mu\cdot \klgs \alpha)\alpha
\leq \ei^{\frac32}\alpha + \kamps\alpha
\tleks
We employ Hardy's inequalities for the harmonic approximation of \(\mu\)
in the terms of the last parentheses:
\[
\olwf \left(\lapgs\mu+\frac12|\klgs \mu|^2\right)\theta^2\alpha
\leq \eaf{\faces}|\mu|\olwf \left|\frac{\lapgs \mu}{\mu}\right|\alpha+
\frac12(\eaf{\faces}|\mu|)^2\olwf\left|\frac{\kl\mu}{\mu}\right|^2
\alpha \leq \eps\left(\koh{3}+\frac12\koh{2}\eaf{\faces}|\mu|
\right)\olwf|\klgs\alpha|^2
\]
provided that
\[
\eaf{\faces}|\mu|\leq \eps \alpha
\]  
We conclude that after elementary manipulations :
\kseks\label{intampest1}
\olwf \vartheta^2|\klgs \alpha|^2 
-\parat{}\left(\olwf \vartheta^2 \alpha^2\right)-\paraf{}\left(
\olwf \vartheta^2h\alpha^2\right)\leq c\ei^{3/2}\olwf 
\vartheta^2\alpha^2
\tleks
Similalry starting from
\begin{equation}\label{pointphest}
-\phi\lapgs \phi-\phi\parat{\phi}=-
 \klgs\mu\cdot\klgs\phi-\frac1\ei\left(
\phi\paraf{\alpha}\paraf{\phi}+2\phi\klgs\alpha\cdot\klgs\phi
\right)
\end{equation}
we get that
\kseks\label{intampest1}
\olwf \vartheta^2|\klgs \phi|^2 
-\parat{}\left(\olwf \vartheta^2 \phi^2\right)-\paraf{}\left(
\olwf \vartheta^2h\phi^2\right)\leq c\ei^{3/2}\olwf 
\vartheta^2\phi
\tleks
We introduce arbitrary test functions in the angular
variables \(\vartheta\) supported in  \(\faces\).
Furthermore in  \eqref{systemap}, multiplying by
\(\theta\) and  using the harmonic approximation \(\harm{\Lambda}\) of
\(\Lambda\) and its  initial form \(\harm{\Lambda}_0\):
\[
|\lapgs\duna|+|\klgs\duna|^2\leq |\lapgs\harm{\duna}|+2|\klgs\harm{\duna}|^2+
\eps
\]
Then for \(\eta=\eaf{\faces}|\duna_0|\) we apply GHI and get that
\kseks
\olwf\left(\left|\lapgs \harm{\Lambda}_0\right|+
2\left|\kl \harm{\Lambda}_0\right|^2+\eps\right)
\vartheta^2  \leq  
\olwf\left(\koh{3}\eta+2\koh{2}\eta^2+\eps\right)|\klgs\vartheta|^2
\tleks
and conclude that for suitable \(\eta>0\)  
\[
\olwf \vartheta^2|\klgs \phi|^2\leq \olwf\left(|\kamps|+\ei\right)\vartheta^2
\]
if we select \(\beta^2_1=C\).  Also we get that:
\[
\olwf \vartheta^2\left(\paraf{\phi}\right)^2\leq C\olwf 
\left(|\kamps|+\ei\right)\vartheta^2
\]
 We will sketch the derivation of estimates derived
from the preceding integral identities. We multiply the 
equation by \(\vartheta^2\) and we obtain that 
\begin{equation}\label{intphasest2}
2\olwf \vartheta^2|\klgs \phi|^2 -
\ei\parat{}\left(\olwf \vartheta^2\phi^2\right)+
2\paraf{}\left(\olwf h\vartheta^2\phi^2\right)+
2\olwf\vartheta^2\left(\paraf{\phi}\right)^2
\leq 
\olwf \vartheta^2\left(2\duna_1\phi^2 + 2\sigma |\klgs\phi|^2\right)
\end{equation}
We select
 \[
\eps_1=\frac14,\qquad \eps_2=\frac14\eaf{\faces}|\alpha|
\]
and obtain  that:
\begin{equation}
2\olwf\vartheta^2|\klgs \phi|^2-\parat{}\left(\olwf\vartheta^2\phi^2\right)+
2\olwf\vartheta^2\left(\paraf{\phi}\right)^2+
2 \paraf{}\olwf h\vartheta^2\phi^2\leq 2\olwf\duna\vartheta^2\phi^2
\end{equation}
and simplifies to the  following inequality:
\kseks\label{phasin}
2\olwf\vartheta^2|\klgs \phi|^2-\ei\parat{}\left(\olwf\vartheta^2\phi^2\right)
+\ei \paraf{}\olwf h\vartheta^2\phi^2 \leq \olwf \duna\vartheta^2\phi^2
\tleks
These are the basic equations that will be used in order to derive the
neccessary estimates. 

\begin{ewfest}
Let   \(\vartheta\) be  a cut off along  the shell
\( \faces\),   then the following holds true:
\[
c_2(r,\faces)){\cal D}^1(\vartheta\phi;\faces)
\leq {\cal D}^0(\vartheta\phi;\faces)\leq
C_2(r,\faces){\cal D}^1(\vartheta\phi;\faces)
\]
\end{ewfest}

\subsubsection{The shell estimates} We observe that the equations \eqref{systemap2},
\eqref{systemap4} are of the form as \eqref{systemap2} and hence are set in
the integral form:
\begin{equation}\label{basic_sh}
\olwf \vartheta^2|\klgs v|^2-\parat{}\left(\olwf\vartheta^2v^2\right)
+\paraf{}\left(\olwf h\vartheta^2v^2\right)\leq \olwf \leda\vartheta^2v^2
\end{equation}
where \(\leda\) is an expression comprising \(\mu,\R,\klgs\R\) etc.
We introduce  the quantities
\begin{equation}\begin{split}
\Pi_0(r) &= \olwf \vartheta^2v^2 \quad\geq 0 \\
\Pi_1(r) &= \olwf  \left(\vartheta\paraf{v}\right)^2\geq 0
\end{split}\end{equation}
and we integrate \eqref{basic_sh} in the interval 
\( \ire=((1-\eps)r,(1+\eps)r)\) 
choosing
\[ 
\fo{(\vartheta)}\cap \{ r/ (r,\theta)\in \shells \} \subset \ire
\] 
and obtain:
\begin{equation}
\begin{split}
\invare\left[ \left(\paraf{\Pi_0(r')}\right)^2
-\paraf{\Pi_0}\ypf{1}\right]dr'+\invare\dirfin{v}{\vartheta}
\Pi_0(r')dr' = \\
=-\invare \Pi_0(r')\ypf{2}(r')dr'
\end{split}
\end{equation}
where 
\begin{subequations}\label{grp}
\begin{align}
\dirfin{v}{\vartheta} &= \olwf \vartheta^2|\klgs v|^2,\label{}\\
\ypf{1}&= \olwf h\vartheta^2v^2,\label{}\\
\ypf{2} &= \olwf \left(\leda\vartheta^2
+\vartheta\parat{\vartheta}+\left(\paraf{\vartheta}\right)^2+
|\klgs \vartheta|^2-\vartheta\lapgs\vartheta\right)v^2
\end{align}
\end{subequations}
We notice first that
\begin{subequations}\label{grp}
\begin{align}
\invare\ypf{1}\paraf{\Pi_0}\leq 
\frac12\invare\left(\paraf{\Pi_0}\right)^2+\frac12 \invare\ypf{1}^2 \label{}\\ 
\ypf{1}\leq C(\shells)\left(\olwf (h\vartheta)^2 \right)^{1/2}\left[
\left(\dirfin{v}{\vartheta}\right)^{1/2}+\left(\olwf
|h|U^2\right)^{1/2}\right]
\end{align}
\end{subequations}
In conclusion we obtain
\[
\invare\left(\paraf{\Pi_0}\right)^2\leq C(\bricks)
\invare\Pi_0^2
\]

\begin{ewfest}
Let   \(\vartheta\) be  a cut off along  the shell
\( \faces\),   then the following holds true:
\[
c_2(r,\faces)){\cal D}^1(\vartheta\phi;\faces)
\leq {\cal D}^0(\vartheta\phi;\faces)\leq
C_2(r,\faces){\cal D}^1(\vartheta\phi;\faces)
\]
\end{ewfest}
\noindent We  will derive an upper bound for \(\psi^{2k},\, \psi=\log\Pi_0\)
then for \(r\in[\delta^{\frac1k},1]\):
\[
|\psi|\leq \frac1\delta\left(|\Pi_0|^\delta+\frac{1}{|\Pi_0|^\delta}\right)
\]  
Specifically we have by Hardy's inequality
\[
\invare \psi^{2kp}\leq C_p(2kr(1+\eps))^p\invare  |\psi|^{(2k-1)p}|\psi'|^p
\leq   \frac{C_p}{\delta^{\frac1\delta}}(2kr(1+\eps))^p\invare
|\psi|^{(2k-1-\frac1\delta)p}|\Pi_0'|^p
\]
for \(0<\delta<1\) that we choose shortly. After H\"older inequality we 
obtain:
\[
\invare \psi^{2kp}\leq 
\frac{C_pC(\bricks)}{\delta^{\frac1\delta}}
(2kr(1+\eps))^p\left(\invare
|\psi|^{(2k-1-\frac1\delta)\frac{2p}{2-p}}\right)
\left(\invare|\Pi_0|^2\right)^{\frac{p}{2}}
\]
Then selecting \(\delta=\frac{1}{k+\eps'-1}, \eps'>0\) we get:
\[
\invare \psi^{2kp} \leq \alpha
\left(\invare|\Pi_0|^2\right)^{\frac{\eps'p}{k}}, \qquad 
\alpha_k=\left(2kr(1+\eps)\right)^{2+\frac{pk}{\eps'}}
\left[C_pC(\bricks)(k+\eps'-1)^{k+\eps'-1}\right]^{\frac{k}{\eps'}}
\]
Then we have that
\[
\invare \psi^{2kp}\leq
\alpha_k^{\frac{2k-\eps'p}{k-\eps'p}}\left(1-\frac{\eps'p}{k}\right)
\left(\frac{2\eps'p}{k}\right)^{\frac{k}{k-\eps'p}}
\]
Summing up we get that:
\[
\sum_{k=0}^\infty
\invare \frac{|\psi|^{2k}}{(2k)!}\leq \left(\sum_{k=0}^\infty
\frac{\alpha_k^{\frac{2k-\eps'p}{k-\eps'p}}}{(2k!)}
\left(1-\frac{\eps'p}{k}\right)
\left(\frac{2\eps'p}{k}\right)^{\frac{k}{k-\eps'p}}\right)^{\frac{1}{p}}
\]
This is majorised  after \(\Gamma\)-function duplication formula  by
\[
\sum_{k=0}^\infty \left[C_pC(\faces)(1+\vareps)r\right]^{k}=c(\faces,r)
\]
Therefore we obtain that:
\[
\sum_{j=0}^\infty\frac{x^{2j}}{(2j)!}\geq \frac12e^x
\] 
and hence 
\[
\invare \Pi_0^2 \leq c(\faces,r)
\]
Following the usual iteration obtained in the appendix we obtain the usual
Harnack inequalities. These are culminated in the following

 \begin{harnshell}
The following estimates hold true:    
\kseksw
\olwf\phi^2\theta^2 & \leq  C_{00}({\shells})\eps\\
\olwf|\klgs \phi|^2\theta^2 & \leq  C_{01}({\shells})\eps\\
\olwf|\klgs^2\phi|^2\theta^2 & \leq C_{02}({\shells})\eps
\tleksw
provided that
\[
\olwf\phi^2,\olwf|\kl\phi|^2,\olwf|\kl^2\phi|^2\geq \eps
\]
\end{harnshell}

\subsection{The lower bound}

The lower bound is obtained by an inductive argument based on the 
reduction to the boundary of a pixel. The one dimensional case 
indicates the method. Assume that we have a function \(\phi\) 
defined in the interval \([0,\mu],\, \phi(0)=\phi(\mu)=0\) satisfying 
\[
c_1\int_0^\mu \phi(x)^2 w(x)dx\leq \mu^2\int_0^\mu \phi'(x)^2w(x)dx\leq 
c_2\int_0^\mu \phi(x)^2 w(x)dx
\]
for a positive weight \(w, 0<c_1<c_2\). Then we prove min-max principle
aalows us to assert that the roots of \(\phi\) in \([0,\mu]\) are at l
east \(c_1\). Then we 
will derive the inequality through the min-max principle, standard
eigenvalue and  isoperimetric inequalities.
The construction is based on the estimates of the preceding section 
that lead to Harnack inequalities for the restriction of the
eigenfunction on the boundary of a geodesic pixel.
We consider first  a domain \(\qwrio_\eps\subset \mangm, n\geq 2\) as 
is described in the appendix on the harmonic approximation.
We drop the \(\eps\) subsctript.

\begin{lelo_in}

Let \(\qwrio \) be a domain with smooth boundary \(\partial\qwrio\) equipped 
with a smooth metric \(\gamma\), induced from  the metric \({\bf g}\). 
Let   \(  \momwb{\psi}\) be a smooth nonnegative function satisfying the estimate for
 \(\gamps{j}=|\klgs^j\psi|^2\):
\kseks\label{sunorekt}
c_{j0}\tau  {\cal D}^0(\gamps{j};\partial \qwrio)\leq
{\cal D}^1(\gamps{j};\partial\qwrio)\leq c_{j1} \tau
{\cal D}^0(\gamps{j};\partial\qwrio),\qquad j=0,1,2
\tleks
Let the zero set  of \(\psi,\,\, \nod{}{\psi}\) be \((n-2)\)-rectifiable. 
Moreover let \(\momb{\phi}\) be such that for \(\tau>0\):
\kseks\label{eswtekt}
c_{30} \tau{\cal D}^0(\phi; \qwrio)\leq 
{\cal D}^1(\phi;\qwrio)\leq c_{31}\tau{\cal D}^0(\phi;\qwrio)
\tleks
and 
\[
{\cal D}^0(\phi-\psi;\partial\qwrio)\leq \eps 
{\cal D}^0(\psi;\partial\qwrio)
\]
Then for \(C_4=C_4\left(\tau, c_{10},c_{11},c_{20},c_{21},c_{30},c_{31}\right)\):
\[
\sharp\{ \partial\qwrio \setminus \nod{}{\psi}\} < C_4
\]
and 
\[
 \Ha{n-1}(\nod{}{\phi})\geq C_0 \tau^{-\frac{n-1}{2}}
\]
\(C_0\) is a numerical constant.

\end{lelo_in}  

\noindent In the appendix we prove that if a smooth function satisfies 
estimates \eqref{sunorekt} then in every connected component of 
\[
\qwrio_\eps=\{ x\in\qwrio/\phi(x)>\eps\}
\]
the following inequalities hold: 
\kseks\label{harnekt}
\eaf{\qwrio_\eps} \phi \leq C_0(\faces) \eps,\qquad 
\eaf{\qwrio_\eps}|\klgs \phi|\leq C_1(\faces)\eps,\qquad 
\eaf{\qwrio_\eps}|\klgs^2\phi|\leq C_2(\faces) \eps
\tleks
for the constants \(C_0,C_1,C_2\) explicitly calculated
Moreover the  harmonic approximation method implies that the function   
\(\kappa=\phi-\harm{\phi}\) satisfies the estimate
\[
 \olww{} |\kl(\zeta \kappa)|^2  \leq C\tau\olww{}\left(\zeta\phi\right)^2   
\]

\paragraph{The tubular neighbourhood of a nodal set} 
The initial form of the harmonic function \(\harm{\phi}\)
denoted by \( \harm{\phi}_0\),  is of degree \(m \aplt \sqrt{\tau} \):
Let  
\[
 \tube_\eps(\nod{}{\harm{\phi}_0})=\{ x\in\qwrio/|\harm{\phi}_0(x)|\leq\eps\} 
\] 
and use  Hardy's inequality \eqref{gh2} we obtain
\[
\oltube{\eps}{\harm{\phi}_0} \phi^2 \leq
C_H\eps^{\frac2m}\oltube{\eps}{\harm{\phi}_0}|\kl \phi|^2 
\leq c_{21}\eps^{\frac2m}\tau \oltube{\eps}{\harm{\phi}_0}(\zeta\phi)^2 \leq
C\tau\eps^{\frac2m}\ogkos{{\tube_\eps(\harm{\phi}_0)}}
\mkf{\tube_\eps(\harm{\phi}_0)}|\phi|
\]
The usual Moser iteration gives us that
\[
\eaf{\tube_\eps}|\kappa|\leq
C\tau\eps^{\frac1m}\left(\frac{1}{\ogkos{\tube_\eps(\harm{\phi_0})}}
\oltube{\eps}{\harm{\phi}_0}
|\kappa|^2\right)^{1/2}
\]
Then according to the usual Harnack estimates we have that if
\(\phi_\eps>0\) near 
y\[
\tube_\eps(\harm{\phi}_0)=\{ x\in\qwrio:|\harm{\phi}_0(x)|\leq \eps\} 
\bigcap \{x\in\qwrio: \phi(x)>\eta\eps\}
\]
we have that:
\[
\eaf{\tube}|\kappa_\eps|\leq C \tau\eps^{\frac1m}\eta\eps  
\]
we conclude then that \( \phi\sim \harm{\phi}\) near \(\nod{}{\harm{\phi}}\). 
The harmonic approximation applied on the slice
allows us to construct the tubular neighbourhoods of nodal sets by the
Lojasiewicz inequality for the approximating function. Specifically, we have
that for suitable choice of \(\eps\) and the tube near the singularities of
multiplicity \(m\) of the nodal set:
\[
|\phi|\leq |\kappa|+|\harm{\phi}|\leq C(\tau\eta\eps^{\frac1m}+1)\eps 
\]
and selecting \(\tau\eps^{\frac1m}=1,\eta=2\)
\[
\nod{}{\phi}\subset \tube_{3\tau^{-m}}(\harm{\phi})
\]

\paragraph{The inductive argument} For \(n=2\) we reduce on a disc and
then we derive estimates for the zero sets using the functions \(\harm{u}\) 
in order to produce test functions for the application of the min-max 
theorem, as it was used in the Courant nodal domain theorem. We recall here
the following lemma from \cite{hs} 

\begin{noco}
There exists \(\eta=\eta_0(n)\in(0,\frac12]\) such that with
\(\eta\in(0,\eta_0)\) if \(w_1,w_2\in C^2(B_2(0)), |w_j|_{C^2}\leq 1\) and
if \(|w_1-w_2|_{C^1}\leq \frac{\eta^5}{2}\) then
\[
{\cal H}^{n-1}(B_{2-\eta}(0)\bigcap \{w_1=0,|\kl w_1|\geq \eta\}) 
\leq (1+c\eta){\cal H}^{n-1}(B_2(0)\bigcap \{w_2=0,|\kl w_2|\geq
\frac{\eta}{2}\}) 
\]
\end{noco}
\paragraph{Estimates on nodal domains and eigenvalues}
We recall the definition of higher order eigenvalues as
\[
\ei_k=\max_{S_{k-1}\subset H^1(\partial\qwrio)}\min_{v\in S_{k-1}^\perp}
\left(\frac{{\cal D}^1(v;\partial\qwrio)}{{\cal D}^0(v;\partial\qwrio)}\right)
\]
The min-max principle for the eigenvalue, \(\ei_k\) suggests also that
\[
\ei_k=\min_{S_k\subset H^1(\partial\qwrio)}
\max_{v\in S_k}
\left(\frac{{\cal D}^1(v;\partial\qwrio)}{{\cal D}^0(v;\partial\qwrio)}\right)
\]
and therefore:
\[
\ei_k\leq  \max_{v\in S_k}
\left(\frac{{\cal D}^1(v;\partial\qwrio)}{{\cal D}^0(v;\partial\qwrio)}\right)
\]
\paragraph{Upper bounds on eigenvalues.} Let us  denote that 
\[
\sharp\{ \partial \faces \setminus \nod{}{\psi}\}=k 
\]
and having selected \(\faces \) containing a geodesic disc of radius 
\(\frac{1}{\sqrt{\tau}}\) then \(k>1\).
We select as  \( S_k= \{\zeta_1,\dots,\zeta_k\} \) where
\(\zeta_j:\qwrio\rightarrow \xw{}, \zeta_j>0 \) defined as follows. 
Let \(\{\cour{i}\}_{j=1}^k\) be the  nodal domains of \(\psi\). 
Set then the tubular neighbourhood 
\[
\cour{i,\eps}=\{ x\in\cour{i}/d(x,\nod{}{\psi})>\eps\}
\]
and 
\[
\nod{j,\eps}{\psi}=\partial \cour{j,\eps}
\]
For this we approximated \(\psi\) harmonically and replaced near its nodes
by \(\widehat{\psi}_0\) so that the Hardt-Simon estimate holds. 
We set as 
\[
\zeta_j|_{\cour{j,\eps}}=\tau_j
\]
and
\[
\zeta_j|_{\cour{i}}=0, j\neq i, \quad \mbox{or} \quad 
\zeta_j|_{\tube_{\theta\eps}\nod{}{\widehat{\psi}_0}\cap \cour{j}}=0
\]
and 
\[
|\klgs \zeta_j|\leq \eta_j
\]
for \(\tau_j,\eta_j\) to be selected. After Sard 's lemma 
\(\nod{j,\eps}{\psi}\) is smooth for suitable \(\eps>0\)
and hence we assume that \(\partial\cour{i}\) 
is also smooth 
We approximate harmonically \(\psi\) and construct 
a smooth partition  of unity, each member supported in a  connected component. 
The Harnack inequalities in the appendix suggest that  
\[
||\harm{\psi}-\psi||_{2,\faces} \leq C(\faces) \eps
\] 
and 
\[
\eaf{\partial\qwrio}|\psi|\leq \eta
\]
We compute that
\[
\ei_k \leq
\frac{\sum_{j=1}^k\tilde{b}_j\eta_j^2}{\sum_{j=1}^k(\mu_j^2\tilde{b}_j+b_j)\tau_j^2}
\]
where 
\[
a_j =\ogkos{\nod{j,\eps}{\psi}},\quad b_j=\ogkos{\cour{j,\eps}},
\quad \tilde{b}_j=
\ogkos{\cour{j,\eps}\setminus \cour{j,\mu\eps}}
\]
Sard's lemma again allows us to choose \(\eps,\mu\) so that 
\[
\tilde{b}_j\sim \mu b_j
\]
Furthermore the isoperimetric inequality suggests:
\[
\tilde{b}_j^{\frac{n-2}{n-1}}\leq
C\left(a_j(1+\eps_j+\eaf{\partial\qwrio}|\klgs\psi|)+b_j\eaf{\partial\qwrio}|h|\right)
\]
Moreover for suitable \(\eps>0\):
\[
\sum_{j=1}^k b_j = \ogkos{\qwrio}-\ogkos{\tube_\eps(\psi)}\geq
(1-\eps)\ogkos{\qwrio}-\left(\ogkos{\nod{}{\psi}}\right)^{\frac{n-1}{n-2}}
\]
Let \(b_1,b_k\) be such that:
\[
b_1\geq \frac1k\ogkos{\qwrio}
\]
and similarly 
\[
b_k\leq \frac1k\ogkos{\qwrio}
\]
We conclude
\[
\ei_k\leq
\frac{k\left(\Ha{n-2}(\nod{}{\psi})\right)^{\frac{n-1}{n-2}}}{\eps\ogkos{\qwrio}}
\]
\paragraph{} The first eigenvalue of the laplacian in 
\(\cour{i,\eps}\) satisfies :
\[
\ei_1(\cour{i,\eps})\leq
\frac{\int_{\cour{i,\eps}}|\klgs \psi|^2}{\int_{\cour{i,\eps}}\psi^2}
\]
Therefore 
\[
\ei_1(\cour{i,\eps})\leq c_{01}\tau
\]
We need now  a lower bound for the first eigenvalue of \(\cour{i,\eps}\):
we select  \(\eps\), so that we avoid dumbbell shape of the nodal domain.
Cheeger estimate of the first eigenvalue  combined with the isoperimetric 
inequality on the  spherical piece suggest that
\[
 2c_{01}\tau\geq \ei_1(\cour{i,\eps}) \geq 
\frac{c}{\ogkos{\cour{i,\eps}}^{\frac{2}{n-1}}}-\mkf{\cour{i,\eps}}|h|
\] 
or 
\[
\ogkos{\cour{i,\eps}}\geq \left(c_{01}\tau+\mkf{\cour{i,\eps}}|h|\right)^{-\frac{n-1}{2}}
\]
Hence since at least one nodal domain should have volume at most 
\[
( c_{01}\tau)^{-\frac{n-1}{2}}\leq \frac1k \ogkos{\faces}
\]
hence we have that 
\[
k\leq \ogkos{\faces}(c_{01}\tau)^{\frac{n-1}{2}}
\]

\paragraph{Lower bounds on eigenvalues} 
The max-min part suggests that
\[
\ei_k\max_{S_{k-1}} \min_{\zeta\in S_{k-1}^\perp}\frac{\int_\faces |\kl
\zeta|^2}{\int_\faces \zeta^2}
\]
Therefore in order to construct a test space \(S_{k-1}^\perp\) 
we have to merge some of  the nodal domains of \(\nod{}{\psi}\). 
We start introducing two numbers for  \(l=1,\dots,k\) and for \(\psi\geq 0\):
\[
D_l^2=\frac{1}{\ei N}\int_{\cour{l}}|\kl \psi|^2,\quad 
W_l^2=\frac{1}{N}\int_{\cour{l}}\psi^2,\quad 
P_l=\frac{1}{P}\int_{\cour{l}}\psi, 
\]
where
\[
N=\int_\faces\psi^2 ,\quad P=\int_\faces \psi
\]
We deduce after dyadic considerations\(\{\cour{\ell}\}_{\ell=1,\dots,k}\)'s that we can select two 
domains, for \(\ell=1,2\) and find some \(j_1\geq 1\):
\[   
\frac{1}{2^{j_1}}\leq |D_1^2-D_2^2|\leq \frac{1}{2^{j_1}}+\frac{1}{2^{j_1+1}}
\]
Similarly for \(j_2,j_3\geq 1\):
\[   
\frac{1}{2^{j_2}}\leq |N_1^2-N_2^2|\leq \frac{1}{2^{j_2}}+\frac{1}{2^{j_2+1}}
\]
and
\[   
\frac{1}{2^{j_3}}\leq |P_1^2-P_2^2|\leq \frac{1}{2^{j_3}}+\frac{1}{2^{j_3+1}}
\]
As \(k\) increases then we select the smallest triple \((j_1,j_2,j_3)\) with
this order.  We consider the space of functions of the form for \(\eta\) a
parameter that we select appropriately and compensates the growth of the
function in the two nodal domains:
\[
\Psi=a_1\left(\frac{1}{P_1}\psi\chi_{\cour{1}}-\frac{\eta}{P_2}\psi\chi_{\cour{2}}
\right)+\sum_{j=3,\dots,k} a_j \psi\chi_{\cour{j}}
\]
We compute that 
\[
\Psi^2=a_1^2\left(\frac{1}{P_1}\psi\chi_{\cour{1}}-
\frac{\eta}{P_2}\psi\chi_{\cour{2}}\right)^2
+\sum_{j=3,\dots,k} a_j^2 \psi^2\chi_{\cour{j}}
\]
and 
\[
|\kl\Psi|^2=a_1^2\left(\frac{1}{P_1}\chi_{\cour{1}}\klgs\psi-
\frac{\eta}{P_2}\chi_{\cour{2}}\klgs\psi\right)^2
+\sum_{j=3,\dots,k} a_j^2 |\klgs\psi|^2\chi_{\cour{j}}
\]
Furtheromore we have 
\[
\int_\faces|\kl\Psi|^2=a_1^2\int_\faces
\left(\frac{1}{P_1}\chi_{\cour{1}}\klgs\psi-\frac{\eta}{P_2}
\chi_{\cour{2}}\klgs\psi\right)^2
+\sum_{j=3,\dots,k} a_j^2 D_j^2
\]
Using Lagrange's identity we write the first term as:
\[
\int_\faces
\left(\frac{1}{P_1}\chi_{\cour{1}}\klgs\psi-\frac{\eta}{P_2}\chi_{\cour{2}}\klgs\psi
\right)^2 
=\left(\frac{1}{P_1P_2}\right)^2
\left[\left(D_1^4+\eta D_2^4\right)\left(P_1^4+P_2^4\right)-
\left(P_1^2D_1^2-P_2^2\eta D_2^2\right)^2\right]^{\frac12}
\]
After Young's inequality since  
\(D_{12}^2=D_1^2+D_2^2,P_{12}=P_1+P_2\) and setting 
\[
a=\frac{P_1}{P_{12}},\qquad b=\frac{P_2}{P_{12}}, 
\qquad 
\upsilon=\frac{\eta D_2}{D_{12}}
\]
we deduce that the last term is minimized by the following expression:
\[
4\ei\left[(4ab-a^2b^2-1)\upsilon^2+
2a^2(1-2ab)\upsilon+\frac{1}{16}-a^4\right]^{\frac12}
\]
After a tedious but otherwise elementary calculation we select 
 \(a\) close to \( 1 \) and \(\eta\) sufficiently
big then we arrange that  the Rayleigh quotient is bounded by \(c\ei\), 
for a suitable constant  \(c>0\).

\begin{itemize}

\item 
\[
\frac34\leq D_1\leq 1, \qquad 0\leq D_2<\frac14 
\]
which implies that 
\[
\frac12\leq D_1-D_2\leq 1
\]
\item  
\[
\frac12\leq D_1\leq \frac34, \qquad \frac14 D_2<\frac12 
\]
which implies
\[
0\leq D_1-D_2\leq \frac12
\]
\end{itemize}
Therefore we have to majorise a function of the form:
\[
Q(a_1,\dots, a_{k-1})=
\frac{\sum_{j=1}^{k-1}D_{j-1}a_j^2}{\sum_{j=1}^{k-1}W_ja_j^2}\geq
\frac{\ei}{k}
\]
We pick \(\eps\sim \frac{1}{\sqrt{\tau}}\) and compute that:
\[
k^2\Ha{n-2}(\nod{}{\psi}\cap \faces) \geq \sqrt{\tau}\ogkos{\faces}
\]
Therefore since \(k^2\leq c\tau^{\frac{n-2}{2}}
\Ha{n-2}(\nod{}{\psi}\cap\faces),\ogkos{\faces}\geq \tau^{\frac{n-1}{2}}\): 
we conclude that
\[
\Ha{n-2}(\nod{}{\psi})\geq C \tau^{-\frac{n-1}{2}}
\]
\paragraph{Modification of the pixel for the singularities of \(\phi\)}
Hopf's strong maximum principle for \(\harm{\phi}\) 
guarantees that \(\nod{}{\harm{\phi}}\) meets transversely the boundary of 
the pixel. The comparison of  nodal sets reduces the problem to the
estimation of the nodal set of \(\harm{\phi}\). 
Therefore using the geodesic spheres  starting near
\(\nod{}{\psi}\) we obtain by the coarea formula, away from 
singularities of \(\harm{\phi}\):
\[
\Ha{n-1}(\nod{}{\phi})= \sum_{j=0}^m\int_{\eps_{j+1}}^{\eps_j} dr 
\int_{\nod{}{\psi}\subset S_r}d\Ha{n-2}(\nod{}{\psi})
\geq c\tau^{-\frac{n-1}{2}}
\] 

\paragraph{The eigenfunctions} The eigenfunctions fullfill the hypothesis 
of the preceding theorem and therefore

\begin{lo_in}
Let \( u : \mangm \rightarrow \xw{}, \lap{g}u=-\ei u, \bricks \subset \mangm\) be a pixel. 
Then 
\[
\Ha{n-1}(\nod{}{u_\ei} \cap \bricks ) \geq C(\bricks)\sqrt{\ei}
\]
\end{lo_in}

\noindent We are placed in the  \(\bricks\)  with boundary spanned by the fronts 
\( \faces_\ell, \ell=1,\dots,m_{i_1\dots i_k}\):
\[
\partial\bricks=\bigcup_{k=1}^{n+1}
\bigcup_{\, \mbox{all} \: i_1,\dots, i_k}
\bigcup_{\ell=1}^{m_{i_1\dots i_k}} \faces_\ell
\]
The boundary of the pixel is not smooth, but we apply the smoothing method
described in the appendix.
Near each front we have the representation of the eigenfunction in the form:
\[
u(r,\dian{\theta})=e^{-\frac{1}{\sqrt{\ei}}\mu(r,\dian{\theta})+
\beta_2\alpha(\dian{\theta})}
\sin(\beta_1\phi(r,\dian{\theta}))
\]
In the preceding paragrpah we obtained the estimates that \(\phi\) satisfeis
with the constant \(c_{11}\sim \ei\). Therefore we have that inside a pixel
of diameter \(\frac{1}{\sqrt{\ei}}\):
\[
\Ha{n-1}(\nod{}{u_\ei} \cap \bricks)\geq C(\bricks)\ei^{-\frac{n-1}{2}}
\]
and since the manifold is splitted in \(\ei^{\frac{n}{2}}\) pixels of trhis
size we have the required estimate.
\subsection{The upper bound}

\cite{do} 
proved that the Hausdorff measure of the nodal set
contained in a pixel \(\bricks\) with  its boundary smoothed out 
is majorised by:
\begin{equation}\label{locdong}
 \Ha{n-1}(\nod{}{u_\ei} \cap \bricks) \leq\frac12\olwb |\kl \log q|+\sqrt{n\ei} \mbox{vol}(\bricks)+
\mbox{vol}(\partial\bricks) 
\end{equation}
where 
\[
 q(u)= |\kl u|^2+\frac{\ei}{n}u^2 
\]
We split the set \(\bricks\) in three parts 
\begin{itemize}
\item The part that is free of nodes:
\[
\bricks_\eps=\{ x\in \bricks/ |u(x)|\geq
\eta\}=\tilde{\bricks}\setminus\tube_\eta(\nod{}{u})
\]
\item The tubular neighbourhood of the nodal set
\(\tube_\eta(\nod{}{u}\cap \bricks)\)  is splitted further as 
\[
\tube_\eta(\nod{}{u}\cap \bricks)=
\mathscr{R}_\eta(\nod{}{u}\cap \bricks)\cup 
\mathscr{S}_\eta(\nod{}{u}\cap \bricks)
\]
\subitem The regular part of the nodal set 
\[
\mathscr{R}_\eta(\nod{}{u}\cap \bricks)=
\{ x\in \bricks/ |u(x)|<\eta, |\kl u(x)|\geq \eta\}
\] 
\subitem and the neighbourhood of the singular set:
\[
\sing_\eta(\nod{}{u}\cap \bricks)=
\{x\in \bricks/ |u(x)|<\eta,|\kl u(x)|< \eta \}
\]
\end{itemize}

The problems in \eqref{locdong} arise in the singular part of the nodal
set.  We will estimate the behaviour of \(|\kl \ln q|\) near
\(\sing_\eps(\bricks)\). We will use induction with respect to the 
multiplicity of the nodal set,  introducing new pixels of 
multiplicity  bounded from below and approximate harmonically the
eignefunction there:  
\[
\tube_\eta(\nod{}{u}\cap\bricks)=
\bigcup_{\ell=0}^{m(\lambda)} \bigcup_{m=1}^{c_\ell} \sing_{\ell,m}
\]
where   \(\sing_{\ell,m}\) is a connected component of mulitiplicity 
\(\ell\)  
\[
\sing_{\eta;\ell,m}=\{ x\in\bricks/ u(x)^2+\cdots+|\kl^{\ell-1}u(x)|\leq \eta,
|\kl^\ell u(x)|\geq \eta\}
\]
Notice that 
\[
\sing_{0,m}=\bricks_\eps,\qquad \sing_{1,m}=\reg_\eta(\nod{}{u}\cap\bricks)
\] 
The harmonic approximation  of \(u_\ei\) in \(\sing_{\ell,m}\) 
is denoted  \(\harm{u}_{\ell,m}\).  We introduce the localization functions 
\(\zeta_{\ell,m}\) with supports
\(\fo(\zeta_{\ell,m})\subset \sing_{\ell,m} \)
We split the integral as
\[
\olwb |\kl \log q|\leq 
\sum_{m=0}^{c_0}\int_{\bricks_\eta}\zeta_{0,m}|\kl\log q|+
\sum_{\ell=0,m=1}^{\lambda,c_\ell}\olwbn{\ell,m} \zeta_{\ell,m}|\kl \log q|
\]
Therefore we set 
\kseksw
I_{\ell,j}[\zeta_m]= \olwbn{\ell,j} \zeta_{\ell,j}|\kl \log q|
\leq C\olwbn{\ell,j}\zeta_{\ell,j} \dong
\tleksw
where we set 
\[
\dong=\frac{|\kl|\kl u|^2|+c|u||\kl u|}{|\kl u|^2+cu^2},\qquad 
\Dong_{\ell,j}=\olwbn{\ell,j} \zeta_{\ell,j}\dong,
\]
We will estimate the terms appearing in the right hand side of the inequality. 
If \(\ell=0\) we have that for \(\hco(\eta,\ei)=\eta^{nt}\ei^{qn+3}\) where 
\(t\) is a parameter to be selected:
\kseksw
\begin{split}
 \eta \leq |u(x)| & \leq \hco(\eta,\ei)\eta\\
|\kl u(x)|& \leq \hco(\eta,\ei)\sqrt{\ei}\eta\\
|\kl^2u(x)| & \leq C\hco(\eta,\ei)\ei\eta
\end{split}
\tleksw
Therefore we have that 
\[
\Dong_{0,j}\leq \ei^{-\frac{n}{2}+qn+\frac32}\eta^{ns}
\]
and we select as 
\[
\eta=\ei^{-ms}, \qquad m=\frac{q}{s}+\frac{1}{ns}
\]
If \(\ell=1\) then we split the tube around the regular \(\mathscr{R}\) of
\(\nod(u)\) in the layers 
\[
\mathscr{R}=\bigcup_{j=1}^\infty\mathscr{R}_j,\qquad 
\mathscr{R}_j= \{ x\in /\theta^j\eta\leq |u(x)|\theta^{j-1}\eta, 
|\kl u|\geq \eta\}
\]
We have that in \(\mathscr{R}_j\)
\kseksw
\begin{split}
 |u(x)| & \leq \hco(\eta,\ei)\theta^j\eta\\
\eta \leq |\kl u(x)|& \leq \hco(\eta,\ei)\sqrt{\ei}\eta\\
|\kl^2u(x)| & \leq C\hco(\eta,\ei)\ei\eta
\end{split}
\tleksw
In this case we get 
\[
\Dong_{0,j}\leq \ei^{-\frac{n}{2}+qn+\frac32}\eta^{ns}
\]
and we select \(\eta\) as before.
If \(\ell>2\): For this we need sharper  estimates 
and we appeal to the harmonic approximation in order to use \L ojasiewicz
inequalities.   We exhaust \(\sing_{\ell,m}\):
\[
\sing_{\ell,m}=\bigcup_{j=0}^\infty \qtube^j_\ell(\theta,\eta),\qquad 
\qtube^j_\ell(\theta,\eta)=
\tube^j_{\ell,m}\setminus\tube^{j+1}_{\ell,m}
\]
for 
\[
\tube^j_{\ell,m}=\sing_{\ell,m}\bigcap \{ x\in\bricks/
\theta^{j+1}\eta\leq |u|\leq \theta^j\eta,\,\,
\theta^{j+1}\eta\leq |\kl u|\leq \theta^j\eta  \}
\]
We introduce in each \(\qtube^j_{\ell,m}(\theta,\eta)\) then
\[
\kappa_{\ell,m;\pm}=u\pm \harm{u}_{\ell,m}
\]
and drop the indices
\[
\kl |\kl u|^2= \kl (\kl \kappa_+\cdot \kl\kappa_-)+\kl|\kl \harm{u}|^2=
(\kl^2\kappa_+)\cdot\kl\kappa_-+(\kl^2\kappa_-)\cdot\kl\kappa_++
\kl|\kl \harm{u}|^2
\]
The harmonic approximation method combined with Harnack inequality 
suggests that in  \(\qtube^j_{\ell,m}(\theta,\eta),\ell>1\) we have 
Bernstein's inequalities 
\[ 
|\kl u|\leq \hco(\eta,\ei)\sqrt{\ei}\eta\theta^j
\qquad 
|\kl^2u|\leq\hco(\eta,\ei)\ei\eta\theta^j
\]
and also 
\[
|\kl\kappa_-|\leq C_1\ei\theta^{2(j+1)}\eta
\] 
Set 
\[
\beta=\hco(\eta,\ei)\eta\theta^j
\]
Furthermore since 
\(|a+b|\geq (1-\varepsilon)\left(a^2 -\frac{1}{\varepsilon}b^2\right) \) 
then :
\[
|u|^2\geq (1-\varepsilon)\left(\harm{u}^2 -\frac1\varepsilon |\kappa_-|^2 
\right)\geq\left(1 -\frac{\sqrt{\ei}\beta}{\varepsilon}\right)\harm{u}^2
\]
and
\[
|\kl u|^2 \geq  (1-\varepsilon)\left(1-\frac{\sqrt{\ei}\beta}{\varepsilon}
\right)
|\kl\harm{u}|^2
\geq(1-\varepsilon)\left(1-\frac{\sqrt{\ei}\beta}{\varepsilon}\right)
|\harm{u}|^{2(1-\nu(\ell))} 
\]
In \(\qtube^j_{\ell,m}\) we have:
\[
u^2\geq \left(1-\frac{\sqrt{\ei}\beta}{\varepsilon}\right)\theta^{2(j+1)}\eta^2,
\qquad 
|\kl u|^2\geq (1-\varepsilon^2)\left(1-\frac{\sqrt{\ei}\beta}{\varepsilon}\right)
\theta^{2(1-\nu(\ell))(j+1)}\eta^{2(1-\nu(\ell))}
\]
The integrand  is estimated according to the estimates derived above:
\[
\zeta_{\ell,m,j}\dong\leq
\zeta_{\ell,m,j}\hco(\eta,\ei)^2\ogkos{\qtube^j_{\ell,m}}
\ei\eta^{2\nu(\ell)}\theta^{2j\nu(\ell)}
\]
We have that due to the multiplicity bound we ahve that
\[
\ogkos{\qtube^j_{\ell,m}}\leq C\ei^{-\frac{n-1}2}\sqrt{\ei}\eta\theta^j, \qquad 
\hco(\eta,\ei)\leq C(\eta^s\ei^q)^n
\]
We set \(\eta=\ei^{-m}\) and compute
\[
m=\frac{4qn+1}{2(2sn+2\nu+1)}+j
\]
 Summing up we get
\[
\sum_{j=0}^\infty\Dong_{\ell,j}\leq C\sum_{j=0}^\infty
\olwbn{\ell,j}\zeta_{\ell,j} 
\leq \frac{C\sqrt{\ei}}{1-\theta^{2\nu}}\ei^{-j}
\]
and therefore we have that
\[
\sum_{j=0}^{C\sqrt{\ei}}\Dong_{j}\leq
\frac{C\sqrt{\ei}}{1-\theta^{2\nu}}
\]We will use repeatedly the following method that we call {\it harmonic
approximation method}. The domains encountered here \(\qwrio\)
 have boundaries with {\it normal crossings} singularities: the 
singular set \(\mathscr{S}(\partial\qwrio)\) is given by the transversal 
intersection of hypersurfaces:  geodesic spheres with local equations 
\(s_1,\dots,s_\ell, \ell=n,n+1\). The piece of the hypersurface 
\(\mathscr{H}_\eps= \{ \dian{x} \in\qwrio 
(s_1\cdots s_\ell+\eps)(\dian{x}=0\}\) near  \(\mathscr{S}\) 
is for suitable \(\eps\)  a smooth hypersurface close to  
\(\mathscr{S}(\partial\qwrio)\). We will consider the domain 
\(\tilde{\qwrio}\) obtained by replacing the singular part
\(\mathscr{S}(\qwrio)\) by  \(\mathscr{H}_\eps\) with repacing the defining
function through cut-offs by the function given there.  
Let \(\harm{F}:\tilde{\qwrio} \rightarrow \xw{}\) be the solution of the 
boundary value problem:
\[
 \lap{} \harm{F} =0,\qquad \harm{F}|_{\partial\tilde{\qwrio}}=F
\]

\paragraph{Harmonic polynomials} We will also approximate the 
harmonic function defined in the pixel \(\harm{F}\) by a sequence 
\(\{F_n\}_{n\in{\bf N}}\) of functions such that 
\begin{equation}\label{iter}
\begin{split}
\Delta_0 \harm{F_0} = 0 ,\qquad
\Delta_0 \harm{F}_n = -\sum_{i,j}{\cal R}^{ij}
\frac{\partial^2 F_{n-1}}{\partial x_i\partial x_j}-g^{ij}\partial_i
\psi\partial_j F_{n-1}
\end{split}
\end{equation}
where
\[
g_{ij}=\delta_{ij}+{\cal R}_{ij} ,\,\, \varrho=\diam{\qwrio},\, 
\psi=\frac12\log(g),\, g=\mbox{det}(g_{ij})
\]
and for \(j=0,1,2\):
\[ 
||\kl^j{\cal R}||\leq C\mu\varrho^{2-j},
\]
Integration by parts after multiplication by  \(\zeta^2F_n\) and
incorporation of the preceding estimates along with Young's inequality  
leads to:
\[
\olww{} \zeta^2|\kl F_n|^2\leq C\mu\varrho^2 \olww{} \zeta^2|\kl F_{n-1}|^2+
C_2\olww{}\left(|\kl\zeta|^2+\zeta^2\right)F_n^2
\] 
and 
\[
\mbox{supp}(|\kl \zeta|)\subset \dact(\qwrio)=
\{x\in \qwrio/d(x,\partial\qwrio)<\eps\}
\]
and
\[
|\kl^j \zeta|\leq \frac{C_j}{\eps^j} 
\]
We select  \(C\mu \rho^2=1\) then
\[
\olww{} \zeta^2|\kl F_n|^2\leq C  \olww{}|\kl(\zeta F_0)|^2
\]
Similarly we have the inequalities:
\[
\olww{}\zeta^2 |\kl^2F_n|^2\leq C\left(\rho^4\olww{}\zeta^2
|\kl^2 F_{n-1}|^2+\rho^2
\olww{}\left(|\lap{} \zeta|+|\kl\zeta|^2\right)|\kl F_{n-1}|^2\right)
\]
and
\[
\olww{}\zeta^2 |\kl^3F_n|^2\leq C^2\left(
\rho^4\olww{}\zeta^2|\kl^3 F_{n_1}|^2+
\rho^2\olww{}\zeta^2|\kl^2F_{n-1}|^2+
\rho^4\olww{}\left(|\lap \zeta|+|\kl \zeta|^2\right)|\kl^2 F_{n-1}|^2\right)
\]
Therefore we have that after iteration:
\[
\olww{} \zeta^2|\kl^2F_n|^2\leq C\olww{}\zeta^2|\kl^2 F_0|^2
\]
and as well as 
\[
\olww{} \zeta^2|\kl^3F_n|^2\leq C\olww{}\zeta^2|\kl^3 F_0|^2
\] 
The Nash-Moser iteration that we describe in the sequel allows us to bound
the sequence in \(C_0^2(\qwrio)\). Rellich lemma allows us to extract a 
sequence that converges in \(H^1(\qwrio)\) and therefore we use it in order
to approximate the initial function by a harmonic polynomial with any accuracy 
we desire. \subsubsection{ The brick localization details} 

Let \(\bricks=\brick{\bjes{\ell}{k}}\) be a brick of size determined by the
parameters \((r(\bricks),\mu(\faces_l)\). We use coordinates 
\( \dian{x}=\dian{\xi}+\dian{c}\) for \(\dian{c}\) denoting the centre 
of the brick.  We construct cut-offs 
Let \(\psi\in C^\infty_0(\overline{\bf R}_+)\) be the following function:
\[
 \psi_\epsilon(t)= \left\{ \begin{array}{lll}
                    1  & \mbox{if} & 0  \leq t\leq 1\\
                    0  & \mbox{if} & t  \geq \frac32
                     \end{array}
                      \right.
\]
Then the following function that localizes in the brick:
\[
\ell_0(\dian{x})=\prod_{i=1}^n 
\psi\left(\frac{\sqrt{4\xi_i^2+\eps_i^2}}{\sqrt{5}\eps_i}\right)
\]
Also we will use the function that localizes in the neighbourhood of the 
zeroes of the  function 
\[
\harm{h}:\bricks\rightarrow \xw{},
\]
\[ 
\nod{}{\harm{h}}=\{\dian{x}\in \bricks:\harm{h}(\dian{x})=0\}
\] 
then the function
\[
\ell_(\dian{x})=\ell_0(\dian{x}_0)\psi\left(\frac{\harm{h}(\dian{x})}{\eps}\right)
\]
localizes in the set:
\[
\nod{\eps}{\harm{h}}=\{\dian{x}\in\bricks/|\harm{h}|\leq \eps\}
\]
We can prove inductively that:
\[
|\kl^j\ell_0|\leq C_j
\]
and 
\[
|\kl^j \ell_{\harm{h}}(\dian{x})|\leq C_j
\sum |\kl^{i_1}\harm{h}|^{p_1}\dots|\kl^{i_m}\harm{h}|^{p_m} 
|\kl^{j_1}\psi|^{q_1}\dots|\kl^{j_l}\psi|^{q_l} 
\]
where we sum over all indices \( i_1p_1+\cdots+i_mp_m+j_1q_1+\cdots+j_lq_l=\ell \)  
and hence that again for indices \( {i_1p_1+\cdots+i_mp_m=\ell,\ell=0,\dots j}\):
\[
|\kl^j\ell_{\harm{h}}(\dian{x})|\leq C_j 
\sum 
\left|\frac{\kl^{i_1}\harm{h}}{\harm{h}}\right|^{p_1}
\dots\left|\frac{\kl^{i_m}\harm{h}}{\harm{h}}\right|^{p_m} 
\]
These will be used succesively in the sequel. 
\subsubsection{\L ojasiewicz, Hardy for functions of the form 
\( \harm{h}\circ N_{\dian{0}}^{-1}\)}

Let \(\harm{h}_0\) be a polynomial function in rectangular coordinates  in 
\(\mpalla{n}{\dian{0},R}\) then we have the following immediate result

\begin{noc}

The function \(\harm{h}=\harm{h}_0\circ N^{-1}\) is a function that satisfies
the following 

\begin{itemize}

\item  The mulitplicity strata \(\Sigma_m(\harm{h})\) 
of the variety \(\nod{}{\harm{h}}\) 
are mapped to \(\Phi(\Sigma_m(\harm{h}_0))=\Sigma(\harm{h})\). 

\item The Lojaziewicz inequalties hold true
\[
|\kl\harm{h}|\geq c_1|\harm{h}|^{1-\ell_1}, \,\,\,\, |\harm{h}(x)|\geq
d(x,\nod{}{\harm{h}}))
\]
  
\end{itemize}
\end{noc}
The first conclusion comes form the chain rule in many variables:
\kseksw
D^\alpha\harm{h}_0= \sum C_{\alpha,\beta,\gamma(j)}
(D^{\beta}\harm{h})(\Phi(x)) (D^{\gamma_1}\Phi)^{e_1}\cdots
(D^{\gamma_l}\Phi)^{e_\ell}
\tleksw
where thje sum extends over all multiindices \(\alpha,\beta,\gamma(j)\in {\bf N}^n,
j=1,\dots,\ell,e_1,\dots,e_\ell\in {\bf N}\) such that:
\[
|\beta|=1, e_1\gamma_1+\cdots+e_\ell\gamma_\ell=|\alpha|-|\beta|
\]
Exchanging the role of \(\harm{h},\harm{h}_0,\Phi,N\) we get the defining equations of
the equimultiple locus.
For the second we  compute:
\kseksw
|\kl f|=|D\Phi_{\dian{x}}(\kl \harm{h}_0)(\Phi(\dian{x})|\geq C 
|\kl \harm{h}_0(\Phi(\dian{x})|
\geq C'|\harm{h}(\Phi(x)|^{1-\ell_1}=C'|f(\dian{x})|^{1-\ell_1}
\tleksw
where we have chosen \(R,r\) so that
\[
|D\Phi(\dian{x})|=|A+B(\dian{x})|\geq |A|-|B(\dian{x})|\geq \frac{|A|}{2}
\]
where \(|B|\leq \frac{|A|}{2}\). Similiarly.
\[
|\harm{h}(\dian{x})|=|\harm{h}_0(\Phi(\dian{x})|\geq c_2 d(\Phi(x),
\nod{\harm{h}_0}{})^{\ell_2}
\] 
now due to the first inequality and the definition of the exponential map
we conclude that:
\[
d(\Phi(\dian{x}),\nod{}{\harm{h}_0})\geq
c'd(\dian{x},\nod{}{\harm{h}}))
\]
A consequence of this is that (GHI)'s hold for such functions.

\subsubsection{Hardy's inequalities}  

Let \(P:\xw{n}\rightarrow\xw{} \) be a homogeneous polynomial of degree \(m\)
and \(\nod{}{P}\) its set of zeroes 
\[
\nod{}{P}=\{ \dian{x}\in\xw{n}/P(\dian{x})=0\}
\]
Moreover let  \(f\in C_0^\infty(\xw{n}\setminus \nod{}{P})\) 
there exist constants, \cite{p}
\( 0<\koh{j}=\frac{4}{(n-2)^2}+O(\eps) n>2\) 

\begin{subequations}\label{gh}
\begin{align}
\int_{\xw{n}} |P|^{-\frac2m} f^2 & \leq \koh{1}\int_{\xw{n}}|\kl f|^2 
\label{gh1}\\
\int_{\xw{n}} \left|\frac{\kl P}{P}\right|^2f^2 & \leq\koh{2}\int_{\xw{n}}|\kl f|^2,
\label{gh2}\\
\int_{\xw{n}} \left|\frac{\Delta P}{P}\right|f^2 & \leq\koh{3}\int_{\xw{n}}|\kl f|^2,
\label{gh3}
\end{align}
\end{subequations}
From the euclidean Hardy' s inequalities we obtain the riemannian
versions by modifying suitably the constants by \(1+\eps\).

\subsubsection{Integration formulas}

Let \(T\) be a tensor field of typr \((p+1,0)\) then we introduce:
\[
A(T)_{i_1\dots i_pk}=T_{i_1\dots i_p;k}-T_{i_1\dots k;i_p}
\]
\[
D(T)_{i_1\dots i_{p-1}}=g^{\ell j}T_{i_1\dots i_{p-1}\ell;j}
\]
The general integration by parts formula reads as

\begin{int_by_parts}

Let \( T\) be a  \((p,0)\) tensor field on the riemannian
manifold \(\mangm\) and \( \phi=|T|^{k-1}\chi,\,\,\chi,\) a smooth 
cut-off  function supported in the domain \(\shells\). Then we have that 
\[
\olws  \phi^2|\kl T|^2 \leq\olws \frac12|A(T)|^2 + |D(T)|^2+
\olws \phi^2\sum_{i=1}^5(\rie*T*T)_i+\olws |\kl \chi|^2 |T|^2
\]
\end{int_by_parts}

\noindent We set \(T_{{\bf i}j}=T_{i_1\dots i_{p-1}j}\) 
and hence we have that 
\[
\phi^2|A(T)|^2=2\phi^2\left(|\kl T|^2 - T_{{\bf i}k;j}T^{{\bf i}j;k}\right)
\]
The last term gives that
\[
\phi^2T_{{\bf i}k;j}T^{{\bf i}j;k}= 
\left(\phi^2T_{{\bf i}k}T^{{\bf i}j;k}\right)_{;j}-
\phi^2 T_{{\bf i}k}T^{{\bf i}j;kj}-2\phi T_{{\bf i}k}\phi_j T^{{\bf i}j;k}
=\mbox{(BT)+(I)+(II)}\]
Furthemore we have that 
\[
(\mbox{I})=\phi^2T_{{\bf i}k}T^{{\bf i}j;kj}=
\phi^2T_{{\bf i}k}\left(D(T)_{{\bf i}}\right)_{;k}
+\phi^2\ric^k_{\,\,s}T^{{\bf i}s}T_{{\bf i}k}+
\phi^2\sum_{\ell=1}^p\rie_s^{\,i_\ell kj}T^{{\bf i}_\ell}T_{{\bf i}k} 
\]
for \({\bf i}_\ell=i_1\cdot i_{\ell-1}si_{\ell+1}\cdot i_p\).
The first term then is written as
\[
 \phi^2T_{{\bf i}k}\left(D(T)_{{\bf i}}\right)_{;k}=(BT)-
\phi^2 |D(T)|^2 - 2\phi T_{{\bf i}j}\phi_j (D(T))^{{\bf i}}
\]
The term (II) is written using that \(\phi=\chi |T|^{k-1}\):
\[
(\mbox{II})=2\phi T_{{\bf i}k}\phi_jA(T)^{{\bf i}jk}+\frac12\kl \log\chi\cdot 
\kl |T|^2+(k-1)\phi^2|\kl |T||^2
\]
In summary we have that:
\[
\olws \phi^2 |\kl T|^2 \leq \olws \phi^2 \left(|D(T)|^2+\frac12|A(T)|^2\right)+
C\olws \left(|\rie|+|\ric|+|\kl \log\chi|^2\right)\phi^2|T|^2
\]

\subsubsection{Functions}

For  a function \(\moms{f},\,\fo{\phi}\subset \shells\) we have that:
\begin{subequations}\label{fint}
\begin{align}
\olws \phi^2|\kl^2 f|^2 &\leq 
C\olws \phi^2\left(|\lap{}f|^2 +|\ric||\kl f|^2\right)\label{fint1}\\
\olws \phi^2|\kl^3f|^2& \leq 
C\olws \phi^2\left( |\kl(\lap{}f)|^2+|\rie| |\kl f|^2\right)\label{fint2}
\end{align}
\end{subequations}

\subsubsection{Curvature}

The diffrential Bianchi identities are written as:
\[
\spcurl(\rie)_{ijkl;m}=\rie_{ijml;k}
\] 
also:
\[
\spapok(\rie)_{ijkl;l}=A(\ric)_{jki}+A(\ric)_{ikj}
\]
We recall here the Bach tensor:
\[
\mbox{B}_{ijk}=A(\ric)_{ijk}+\frac{1}{(n-1)(n-2)}\left[
g_{ij}\kl_k\R-g_{ik}\kl_j\R\right]
\]
and the contracted identities are writen as
\[
\spapok(\ric)_{i}=\klgs_i \R
\]
Therefore we have that
\[
\olws \phi^2|\kl \rie|^2 = 4\olws|B|^2+\frac{2}{(n-1)(n-2)^2}|\kl \R|^2+\sum_{i=1}^5 (\rie*\rie)_i
\]
\paragraph{The iterative method}

Let \( \momw{g,\chi}\) be smooth functions  and \(\moms{\harm{h}}\) a polynomial 
weight function of degree \(m\) then set 
\[
\qwrio_j=\{ x\in \qwrio/ \quad 
\theta(1-\theta^j)\frac{\eta}{2}\leq |\harm{h}(x)| \leq
(1-\theta+\theta^j)\eta\}
\]
and \( \fo{\chi}_j\subset \qwrio_j\). This for instance is given for 
\[
\chi_j(x)=\ell\left(\frac{\harm{h}(x)}{(1-\theta+\theta^j)\eta}\right)
\ell\left(\frac{\theta(1-\theta^j)\eta}{\harm{h}(x)}\right)
\]
We suppose that the smooth function \(g\) satisfies the inequality, for
positive constants \(\gamma>1, e=2,4\) and any smooth cut-off \(\chi\):
\[
\olww{}  \chi^2|\kl g|^2 \leq \gamma\olww{} \chi^2 |g|^e
\]
Then Sobolev inequality  suggests for 
\(s=\frac{np}{n-p},1< p<2, k,\ell,d>1\):
\kseks\label{init}
\left(\olww{j} \left(\chi^d |\harm{h}|^{-\ell} g^ k\right)^s\right)^{p/s}
\leq C_0k \olww{j-1}
\left(|\harm{h}|^{-\ell}\chi^d g^{k-1}\right)^p|\kl g|^p +
\ell\olww{j-1}\left|\frac{\kl \harm{h}}{\harm{h}}\right|^p
\left(\chi^d |\harm{h}|^{-\ell}g^k\right)^p+\\
+p\olww{j-1}\chi^{(d-1)p}|\kl \chi|^pg^{kp}|\harm{h}|^{-p\ell}
\tleks
The first term then gives for \(r=\frac{2p}{2-p}\)
\begin{equation}
\begin{split}
\olww{j}  
\left(|\harm{h}|^{-\ell}\chi^d g^{k-1}\right)^p|\kl g|^p& \leq
\left(\olww{j-1} \chi^{dp}
\left(|\harm{h}|^{-\ell}g^{k-1}\right)^r\right)^{\frac{p}{r}}
\left(\olww{j-1} \chi^{dp} |\kl g|^2\right)^{p/2}\leq \\
 & \leq C_1((1-\theta+\theta^j)\eta)^{p\ell}\gamma^{\frac{p}{2}}
\ogkos{\qwrio_j}^{\frac{(k-1)r-e}{(k-1)r}}
\left(\olww{j-1} \chi^{\frac{2rdp}{e}}|\harm{h}|^{-r\ell} g^{(k-1)r}
\right)^{\frac{p(k-1+\frac{e}{2})}{r(k-1)}}
\end{split}
\end{equation}
since we have  
\kseksw
\left(\olww{j} \chi^{2dp} |\kl g|^2\right)^{p/2} \leq
 \gamma^{\frac{p}{2}}\left((1-\theta+\theta^j)\eta\right)^{\frac{pe\ell}{k-1}}
\ogkos{\qwrio_j}^{\frac{(k-1)r-e}{(k-1)r}}
\left(\olww{j-1} \chi^{\frac{2rdp}{e}}|\harm{h}|^{-r\ell} g^{(k-1)r}
\right)^{\frac{pe}{2r(k-1)}}
\tleksw
 The middle term after  application of \eqref{gh2} gives that:
\kseksw
\begin{split}
\olww{j}\left|\frac{\kl \harm{h}}{\harm{h}}\right|^p
\left(\chi^d |\harm{h}|^{-\ell}g^k\right)^p& \leq
\left(\olww{j-1}
\chi^{dp}\left(h^{-\ell}g^{(k-1)}\right)^r\right)^{\frac{p}{r}}
\left(\olww{j-1}\left|\frac{\kl \harm{h}}{\harm{h}}\right|^2\chi^{dp}g^2
\right)\\
& \leq C_1 ((1-\theta+\theta^j)\eta)^{p\ell}\gamma^{\frac{p}{2}}
\ogkos{\qwrio_j}^{\frac{(k-1)r-e}{(k-1)r}}
\left(\olww{j-1} \chi^{\frac{2rdp}{e}}|\harm{h}|^{-r\ell} g^{(k-1)r}
\right)^{\frac{p(k-1+\frac{e}{2})}{r(k-1)}}
\end{split}
\tleksw
In an analogous way the  last term leads to:
\[
\olww{j} |\kl \chi|^p(\harm{h}^{-\ell}g)^k\leq  
C_1((1-\theta+\theta^j)\eta)^{p\ell}\gamma^{\frac{p}{2}}
\ogkos{\qwrio_j}^{\frac{(k-1)r-e}{(k-1)r}}
\left(\olww{j-1} \chi^{\frac{2rdp}{e}}|\harm{h}|^{-r\ell} g^{(k-1)r}
\right)^{\frac{p(k-1+\frac{e}{2})}{r(k-1)}}
\]
In summary we arrive at
\[
\left(\olww{j} 
\left(\chi^d|\harm{h}|^{-\ell}g^k\right)^s\right)^{\frac{1}{ks}}\leq
C_3k^{\frac1k}(\eta(1-\theta+\theta^j))^{\frac{p\ell}{k}}\gamma^{\frac{p}{2k}}
\ogkos{\qwrio_j}^{\frac{(k-1)r-e}{(k-1)kr}}
\left(\olww{j-1} \chi^{\frac{2}{e}dr}|\harm{h}|^{-r\ell} g^{(k-1)r}
\right)^{\frac{(k-1+\frac{e}{2})}{k}\cdot\frac{1}{r(k-1)}}
\]
Notice that 
\[ 
r=qs, \quad q(n,p)=\frac{1-\frac{p}{n}}{1-\frac{p}{2}}, \qquad 
1-\frac{2}{n}\leq q(n,p)\leq 2-\frac{2}{n}
\]
We conclude with the basic inequality that we will iterate
\[
\left(\frac{1}{v_j}\olww{j} G^{ks}\right)^{\frac{1}{ks}}\leq C_j
\left(\frac{1}{v_{j-1}}\olww{j-1} G^{(k-1)r}\right)^{\frac{1}{(k-1)r}}
\]
where 
\[
r=\frac{s}{a}, \qquad a=\frac{n}{n-1},  \qquad v_j=\ogkos{\qwrio_j}
\]
\[
C_j=C_3\left[k\eta^{\frac{pks+1}{ks}\ell}
\left(\theta(1-\theta^{j-1})\right)^{-\frac{k-1+\frac{e}{2}}{k}\cdot\frac{\ell}{k-1}}
\gamma^{\frac{p}{2k}}
v_j^{\frac{(k-1)r-e}{(k-1)r}}\beta_j\right]^{\frac1k}, 
\]
\[
\beta_j=\frac{v_j^{\frac{1}{r(k_j-1)}}}{v_{j+1}^{\frac{1}{k_{j+1}s}}}
\]
and we replace \(\beta_j\) by the upper bound 
\[
\beta_j\leq \frac{ \eta^{\frac{1}{a(n-2)}}  }{\theta^{\frac1a}}=\sigma
\] 
In order to bring this to the standard iteration form we do the following:
\[
k_j=\frac{a^{j+1}}{s}+1
\]
Finally we arrive at the basic iteration inequality:
\[
I_{j+1}=\left(\frac{1}{v_{j+1}}\olww{j+1}G^{a^{j+1}}\right)^{\frac{1}{a^{j+1}}}\leq
C_j
\left(\frac{1}{v_j}\olww{j}G^{a^j}\right)^{\frac{1}{a^j}},
\]

Then we have the iteration inequality 
\kseks\label{anadr}
I_{j+1} \leq C_j I_j
\tleks
The iteration leads to the inequality 
\[
\eaf{\qwrio_\infty}|G|\leq D \left(\frac{1}{\ogkos{\qwrio_0}}
\olww{0}G^p\right)^{\frac1p},\qquad
D=\lim_{j\rightarrow\infty}\left(C_1\dots C_j\right)
\]
We  select the local density parameter as
\(\theta\sim\frac{1}{\gamma^t},t>0\).
The constant is estimated through  elementary inequalities of the form:
\[
\frac{1}{\gamma^t a^2-1}\leq 
-\sum_{j=0}^\infty \frac{1}{a^{2j}}\log\left(1-\frac{1}{\gamma^{tj}}\right)
\leq \frac{\gamma^t}{\gamma^t-1}\frac{1}{\gamma^ta^2-1}
\]
We arrive finally at 
\[
c=\eta^{\ell pn+\frac{a+1}{a(n-2)}}\gamma^{\frac{pn(t+1)}{2t}+\frac{t-1}{n-2}}
\]
We will denote the constant in the form:
\[
\hco(\eta,\gamma)=\left(\eta^s\gamma^q\right)^n ,
\qquad s=\ell p+\frac{a+1}{na(n-2)},
\qquad
q=\frac{p(t+1)}{2t}+\frac{t}{n(n-2)}
\]

\paragraph{Back to harmonic approximation} The harmonic
approximation method suggests that
\[
\eaf{\qwrio_\infty}|\kappa|\leq \hco(\eta,\ei)
\left(\frac{1}{\ogkos{\qwrio_0}}\olww{0}u^2\right)^{1/2}
\]
Similarly we have the higher order inequalites
\[
\eaf{\qwrio_\infty}|\kl\kappa|\leq
\hco(\eta,\ei^2)
\left(1+\hco(\eta,\mu)\ei\frac{1}{\ogkos{\qwrio_0}}\olww{0}|\ric|\right)^{1/2}
\left(\frac{1}{\ogkos{\qwrio_0}}\olww{0}u^2\right)^{1/2}
\]
Let \(\harm{u}>0\) be a harmonic function. Then set for \(j=0,1\)
\[
h_j=\left|\frac{\kl^{j+1}\harm{u}}{|\kl^j\harm{u}|}\right|^2,\qquad
H_0=\left|\frac{\kl|\kl \harm{u}|^2}{|\kl \harm{u}|^2}\right|
\]
and compute
\[
\lap{g}h_0\geq -(|\ric|+2)h_0^2+\frac23h_1h_0\geq -(|\ric|+2)h_0^2 
\]
or that
\[
\olww{} \zeta^2|\kl h_0|^2\leq  \olww{} (|\ric|+2)\zeta^2h_0^3
\]
We note that the integral is of the form 
\[
\olww{}\phi^2h_0^3\leq (1+\eps^2)\koh{2}\olww{}\phi^2|\kl h_0|^2+
\frac{1}{\eps^2} |\kl \phi|^2
\]
majorised after application of Hardy's inequality. 
We compute that
\[
|\kl h_0|^2\leq (1+\eps)^2(H_0h_0^2+\frac{1}{\eps^2}h_0) 
\]
then we find
\[
\olww{}\phi^2|\kl h_0|^2\leq  (1+\eps)^2\olww{} 
\phi^2(H_0h_0^2+\frac{1}{\eps^2}h_0) \leq
(1+\eps^2)\koh{2}\olww{}|\kl\phi|^2\left(\frac{1}{\eps^2}+h_0^2\right)+ 
\phi^2|\kl h_0|^2
\]
If \(n>1,\eps'>0\) then select  
\[
\frac{n^2-2n-\eps'^2}{n^2-2n+8}=\eps^2<\frac{n^2-2n}{n^2-2n+8}
\]
and conclude that 
\[
\olww{}\phi^2|\kl h_0|^2\leq c_n(\eps')
\olww{}|\kl\phi|^2\left(\frac{1}{\eps^2}+h_0^2\right), \quad
 c_n(\eps')=\frac{8-\eps'^2}{n^2-2n-\eps'^2} 
\]
and hence 
\[
\olww{} \phi^2h_0^3 \leq c_n(\eps')\olww{} \phi_1^2(1+\eps'h_0^2),\quad
\phi_1=|\kl\phi|
\]
Repeating the procedure we get 
\[
\olww{}\phi_1^2h_0^2\leq c_n(\eps')\olww{} 
(1+\eps^2h_0)|\phi_2|^2
\]
and 
\[
\olww{}\phi_2^2h_0\leq c_n(\eps')\olww{} 
(1+\eps^2)|\phi_3|^2,\qquad \phi_3=|\kl\phi_2|
\]
We apply this formula for \(\phi=(|\ric|+2)^{1/2}\zeta\) that
\[
\olww{}\zeta^2|\kl h_0|^2\leq \olww{} 
\sum_{i=1,j\leq i}^3c_n(\eps')^i|\kl^{i-j}\ric||\kl^j\zeta|^2
\]
and we obtain
\[
\eaf{\qwrio}h_0\leq \hco(t,\gamma^{-3})
\olwwr{\qwrio}|\R|, \qquad t=\max(t_1(\qwrio),t_2(\qwrio),t_3(\qwrio) 
\]
where \(|\kl\zeta|^j\leq c\gamma^{-j}\)\paragraph{Growth of  a function near its  nodal set } 

We assume that
\[
\olww{}\zeta^2|\kl u|^2\leq \tau\olww{}\zeta^2u^2
\]   
Let 
\[
\qwrio_j=\{x\in \bricks/ (1-\theta^j)\theta\eta\leq 
|\harm{u}(x)|\leq (1-\theta+\theta^j)\eta \}\bigcap
B_{C,\frac{1}{\sqrt{\tau}}}
\]
and
\[
\dact(\qwrio_j)=
\qwrio_j\setminus
\{x\in \bricks/ (1+\theta^{j+1})\theta^2\leq |\harm{u}(x)|\leq
\theta(1-\theta^j)\eta\}\bigcap B_{C,\frac{1}{\sqrt{\tau}}}
\]
The cut-off  function \(\zeta\)  satisfies the following estimate:
\[
|\kl^\ell\zeta|\leq \frac{c_\ell}{\theta^\ell}
\]
We apply Hardy's inequality
\[
\olww{j}\zeta^2u^2=\olww{j}\harm{u}^{\frac2m-\frac2m}(\zeta u)^2\leq 
\koh{1} \left(\eaf{\qwrio_j}|\harm{u}|\right)^{\frac2m}\olwwr{j}|\kl(\zeta u)|^2 \leq
\]
\[
\leq
C \left(\eaf{\qwrio_j}|\harm{u}|\right)^{\frac2m}(1+\eps)\left[
\frac1\eps\olwwr{\qwrio{j}}|\kl\zeta|^2u^2+
\tau\olwwr{\qwrio{j}}(u\zeta)^2\right]
\]
Therefore we have that close to \(\tube_\eta\left(\nod{}{\harm{u}}\right)\), for 
\(\eta=(2\koh{1}\tau(1+\eps))^{-\frac{m}{2}}\). We have that
\[
\olww{j}\zeta^2u^2\leq
\frac{1}{\tau^{\frac{m}{2}}\eps}\olwwr{\qwrio_j}|\kl \zeta|^2u^2\leq 
\frac{c_1}{\tau^{\frac{m}{2}}\eps\theta^2}\olwwr{\qwrio_j}u^2
\] 
We have that:
\[
\olwwr{\qwrio_j}u^2\leq 2\olwwr{\qwrio_j}\harm{u}^2+\kappa^2
\]
and the second integral is estimated again as
\[
\olwwr{\qwrio_j}\kappa^2\leq
\koh{1}(\theta\eta)^{\frac2m}\tau\olwwr{\qwrio_j} u^2
\]
and hence
\[
\olwwr{\qwrio_j}u^2\leq \frac{4}{2-\theta^{\frac2m}}
\olwwr{\qwrio_j}\harm{u}^2
\]
The coarea formula suggests then after \L ojasiewicz inequality that
\[
\olwwr{\qwrio_j}\zeta^2\harm{u}^2= 
\int_{(\theta-\theta^j)\eta}^{(\theta-\theta^j)(1+\xi)\eta}
d\mu\int_{\{ \harm{u}=\mu\}} \frac{\harm{u}^2d\sigma_\mu}{|\kl \harm{u}|}
+\int_{(1-\theta+\theta^j)\xi\eta}^{(1-\theta+\theta^j)\eta}
d\mu\int_{\{ \harm{u}=\mu\}} \frac{\harm{u}^2d\sigma_\mu}{|\kl \harm{u}|}
\leq 
\]
\[
\int_{(\theta-\theta^j)\eta}^{(\theta-\theta^j)(1+\xi)\eta}
\mu^{\nu+1}\alpha(\mu)d\mu 
+\int_{(1-\theta+\theta^j)\xi\eta}^{(1-\theta+\theta^j)\eta}
\mu^{\nu+1}\alpha(\mu)d\mu
\]
Inside a ball of radius \(\tau\) applying Crofton formula if the
multiplicity of \(\harm{u}\) is \(m\)
\[
\alpha(\mu)=\int_{\{ \harm{u}=\mu\}}d\sigma_\mu\leq c\tau^{-\frac{n-1}{2}}m
\]
Therefore we find that for \(\sigma_j=\theta-\theta^j\):
\[
\olwwr{\qwrio_j}\zeta^2\harm{u}^2\leq
\tau^{-\frac{n-1}{2}}\eta^{\nu+2}
\left(\sigma_j^{\nu+1}\xi^{\nu+2}+
\frac{(1-\sigma_j)^{\nu+2}(1-3\sigma_j)}{1-2\sigma_j} \right)
\]
Finally we have that
\[
\olww{j}\zeta^2u^2\leq C\tau^{-\frac{n-1}{2}}\eta^{\nu+2}
\]
Moreover we recall the folowing identity  from \cite{p}:
\[
\eta^3\eta\frac{d}{d\eta}\eta^{-3}\olww{}\zeta^2\harm{u}^2=
-\olww{}\frac{\kl Q\cdot\kl\harm{u}}{Q^2}\zeta^2\harm{u}^3
\]
for \(Q=|\kl \harm{u}|^2\geq c|\harm{u}|^{2(1-\nu)}\). This  leads to 
\[
\eta^3\eta\frac{d}{d\eta}\eta^{-3}\olww{}\zeta^2\harm{u}^2\leq
\olww{}\left|\frac{\kl Q}{Q}\right|\zeta^2\harm{u}^{\nu+2}\leq 
\olww{}\left|\frac{\kl Q}{Q}\right|^2\zeta^2\harm{u}^{\nu+2}
\]
We apply Hardy's inequality and get 
\[
\eta I'(\eta)\leq
C\left(1+\tau^{-\frac{n+1}{2}}\eta^{\nu}\right)\tau\eta^\nu I(\eta)
\]
for 
\[
I(\eta)=\eta^{-3}\olww{}\zeta^2\harm{u}^2
\]
Finally we get that for \(\eta>\eta_0\):
\[
I(\eta)\leq C
e^{c_\nu\left(1+\tau^{-\frac{n+1}{2}}\eta^{\nu+1}\right)\tau\eta^\nu}
I(\eta_0)
\]
 and 
\[
\olww{}\zeta^2\harm{u}^2 \leq C
e^{c_\nu\left(1+\tau^{-\frac{n+1}{2}}\eta^{\nu+1}\right)\tau\eta^\nu}
\left(\frac{\eta}{\eta_0}\right)^3 \olww{0}\zeta^2\harm{u}^2
\]

\paragraph{Morrey estimates}

Let  \(\eps<1,\,\, 0<\gamma<1\,\, \mbox{or}\,\, \gamma<0,\, p<2\):
\[
u_\eps=\sqrt{u^2+\vareps^2},\qquad 
\psi_\vareps=\log u_\vareps, \qquad w=u_\eps^\gamma 
\]
and for  \(\zeta,\fo{(\zeta)}\subset \qwrio\):
\kseks\label{assum}
\olww{} \zeta^2|\kl u_\eps|^2\leq \tau\olww{} \zeta^2 u_\eps^2
\tleks
Then for \(q=\frac2\gamma\):
\begin{subequations}\label{moiq}
\begin{align}
\olww{} |\kl w|^p\zeta^p &\leq C_1(\tau)\olww{} |\kl \zeta|^2w^q\label{moiq1}\\
\olww{} |\kl \psi_\vareps|^2\zeta^2 &\leq C_2(\tau)\olww{} |\kl
\zeta|^2+\zeta^2\label{moiq2}
\end{align}
\end{subequations}
The inequality \eqref{moiq1}  follows after selection for 
\(\zeta\) as \(u_\eps^{\gamma-1}\) gives since
\[
|\kl w|=\gamma w^{1+\frac1\gamma}|\kl u_\eps|
\]
that 
\[
\olww{} \zeta^2|\kl w|^2\leq \frac{C_0}{\gamma^2}\olww{} \zeta^2w^2
\]
\noindent The inequality \eqref{moiq2} requires the additional assumption
for \(\tau>0\):
\[
\olww{}\zeta^2|\kl^2u|^2\leq \tau^2 \olww{} \zeta^2 u^2
\]
We start selecting values \(v_1,\dots, v_m>0\) and assume that
\[
u=v_j+h_j
\]
making the following choice:
\[
|h_j|\leq \eps |v_j|
\]
then 
\[
(v_j+h_j)^2 \geq (1-\eps)\left(v_j^2-\frac{h_j^2}{\eps}\right)\geq
\left(1-\eps\right)^2h_j^2
\]

We approximate harmonically \(h_j\) in suitable bricks selected 
so that we use the initial form of \(\harm{h}_j\). Hence we have that
for \(\psi=\log(u_\eps),\tilde{\psi}=\log(h)\)
\[
\olwb{} |\kl\psi|^2\zeta^2\leq c\olwb{}|\kl \tilde{\psi}|^2\zeta^2
\]
Therefore  let  \(\harm{h}\) be the harmonic approximation of \(h\) in 
\(\qwrio\)  and we set:
\[
h=\harm{h}+\kappa
\]
The standard harmonic approximation method estimates  combined 
with partial integration leads us to 
\[
\olww{}\zeta^2 |\kl^2\kappa|\leq \hco(\eta,\tau^2)\olww{} \zeta^2u^2
\]
The estimate of the preceding paragraph
\[
|\kl\kappa|\leq \eaf{\qwrio_0}|\kl\kappa|\leq c
\left(\olww{} \zeta^2u^2\right)^{1/2},\qquad
c=\frac{\hco(\eta,\tau)+\hco(\eta,\mu)\tau ||\ric||_{1}}{\ogkos{\qwrio{0}}} 
\]
We compute for \(\eps<1\):
\[
u_\eps^2=\harm{u}^2+2\kappa\harm{u}+\kappa^2+\eps^2\geq
(1-\eps^2)\harm{u}^2+(1-\frac{1}{\eps^2})\kappa^2+\eps^2
\]
Then we select \(\eps\) such that
\[
\kappa^2(1-\frac{1}{\eps^2})+\eps^2>\frac12\eps^2
\]
and therefore we get that
\[
|\kappa|\leq \frac{\eps^2}{\sqrt{2(1-\eps^2)}}
\]
Let then \(m\) denote  the highest multiplicity of 
\(\harm{u}\). We apply the preceding estimates to conclude that for 
\(\eps'=\frac{\eps^2}{2(1-\eps^2)}\):
\[
\olww{} \zeta^2\frac{|\kl u_\eps|^2}{u_\eps^2}\leq 
\frac{1}{2(1-\eps'^2)}\olww{}
\left|\frac{\kl\harm{u}_{\eps'}^2}{\harm{u}_{\eps'}^2}\right|\zeta^2+
c \olww{} \zeta^2u^2\olww{}\frac{\zeta^2}{\harm{u}_{\eps'}^2}
\]
Now we use  the estimate of the preceding paragraph 
\[
||\zeta u||_{2,\qwrio}\leq C\mu\tau^{-\frac{n-1}{2}}\eta^{\nu+2}
\]
and obtain that 
\[
\olww{}\zeta^2|\kl \psi_\eps|^2\leq C_\eps\olww{}|\kl\zeta|^2+
c\tau^{-\frac{n-1}{2}}\eta^{\nu+2}\olww{}\frac{\zeta^2}{u_{\eps'}^2}
\]
We set  \(e_k(m)=1-\frac{k}{m}\)  and then 
\[
\olww{} \frac{\zeta^2}{\harm{u}^2}= 
\olww{} \frac{(\harm{u}^{-e_1}\zeta)^2}{\harm{u}^{\frac2m}}
\leq
\koh{1}\left[e_1^2\olww{} \left|\frac{\kl \harm{u}}{\harm{u}}\right|^2
\left(\harm{u}^{-e_1}\zeta\right)^2 +\olww{}\harm{u}^{-2e_1}|\kl\zeta|^2\right]
\]
The constant in Hardy's inequality is 
 \(\koh{1,2}\sim\frac{4}{(d-1)^2}, d\geq 3\):
\[
\olww{}|\kl\harm{\psi}|^2\harm{u}^{-2e_1}\zeta^2\leq(1+\eps)\koh{1} 
\left[\frac1\eps\olww{}\harm{u}^{-2e_1}|\kl\zeta|^2+e_1^2\olww{}
|\kl\harm{\psi}|^2\harm{u}^{-2e_1}\zeta^2\right]
\]
We select then \(\eps\) so that :
\[
\eps\leq \frac{(d-3)m+2}{(d+1)m-2}
\]
Hence we have 
\kseksw
\olww{}\zeta^2|\kl \harm{\psi}|^2\leq \frac{(1+\eps)\koh{1}}{\eps}
\olww{}\harm{u}^{-2e}|\kl\zeta|^2
\tleksw
and we conclude that:
\kseks\label{mobast}
\olww{}\zeta^2|\kl\psi|^2\leq \frac{(1+\eps)\koh{1}}{\eps}\olww{}
\harm{u}^{-2e}|\kl\zeta|^2
\tleks
Near the zeros of \(u,\harm{u}\), i.e. selecting \(\eta\) suitably small
then we retrieve the same inequality.
Further treatment is required when we introduce in \eqref{mobast} 
\[
\zeta = \psi^\ell\vartheta 
\]
and we get:
\[
\olww{}\vartheta^2
\psi^{2\ell}|\kl\psi|^2\leq C\ell^2
\olww{} \harm{u}^{\frac2m}\left(
\vartheta^2\psi^{2(\ell-2)}|\kl \psi|^2+\psi^{2(\ell-1)}|\kl \vartheta|^2\right)
\]
through the elementary inequality for \(x\in(\delta^{\frac1k},1]\):
\[
|\log x|\leq
\frac2\delta\left(|x|^{\frac{\delta}{2}}+|x|^{-\frac{2}{\delta}}\leq
\frac4\delta^{-\frac{2}{k\delta}-1}\right)
\]
We obtain for \(\varrho<1,\delta=(e\varrho\ell)^{-1}\):
\kseks\label{finid}
\olww{}\vartheta^2\psi^{2\ell}|\kl \psi|^2\leq
2C\ell^{2(e\varrho\ell+1)}\olww{} \psi^{2\ell(1-e\varrho))}
\left[\vartheta^2|\kl\psi|^2+\psi^{2e\varrho\ell}|\kl\vartheta|^2\right]
\tleks

\paragraph{The iteration for the lower bound.} We follow the method of
\cite{lspde}. We select as  \(\zeta\)  
\[
v^{2(\ell-1)}\vartheta^{2(a\ell-b)}, \qquad v=\psi_\eps-\sigma,\qquad a\geq 1+b
\]
obtaining that:
\[
\olww{} \vartheta^{2(a\ell-b)}v^{2(\ell-1)}|\kl v|^2\leq 
\ell^{2(e\varrho\ell+1)}\olww{} v^{2\ell-4}|\kl v|^2
\eta^{2(a\ell-b)}+v^{2(\ell-1)}\vartheta^{2(a\ell-b-1)}
\]
We get that
\[
\olww{} |\kl (|v|^\ell\eta^{a\ell-b})|^2\leq \ell^{2\ell(e\varrho+3)}
\ogkos{\qwrio}+\ell^{2\ell(e\varrho+1)}\olww{} v^{2\ell}\vartheta^{2(a\ell-b-1)}
\]
Therefore we have that for \(a=kb, \qquad b=\frac{m}{k-m},\qquad 
m=\frac{\ell}{\ell-2}\)
\[
\left(\frac{1}{\ogkos{\qwrio}}\olww{}(|v|\vartheta)^{2\ell k}\right)^{\frac1k}
\leq \ell^{2\ell(e\varrho+3)}+\ell^{2\ell(e\varrho+1)}
\left(\frac{1}{\ogkos{\qwrio}}\olww{}(|v|\vartheta^a)^{2m\ell}\right)^{\frac1m} 
\]
Hence we have that 
\[
\left(\frac{1}{\ogkos{\qwrio}}\olww{}(|v|\vartheta)^{2\ell k}\right)^{\frac{1}{\ell k}}
\leq \ell^{2(e\varrho+1)}+\ell^{2(e\varrho+1)}
\left(\frac{1}{\ogkos{\qwrio}}\olww{}(|v|\vartheta)^{2m\ell}\right)^{\frac1m} 
\]
Selecting a sequence \(\{\ell_j\}\):
\[
\ell_j=k^j \qquad 
I_j =\left(\frac{1}{v_j}\olww{}(|v|\eta)^{2k^j}\right)^{1/k^j}
\] 
 Hence we obtain 
\[
I_{j+1}\leq k^{j(e\varrho+3)}+k^{j(e\varrho+1)} I_j
\]
or that
\[
I_{j+1}\leq k^{(e\varrho+3)j}(1+I_j)
\]
Examining separately the two cases: for some \(j_0\):
\[
I_{j_0}\leq 1
\] 
and the complementary case we arrive at the conclusion
\[
I_j\leq (2k^{j_0})^{j-j_0}
\]
Therefore following the reasoning in \cite{lspde} we conclude that 
\[
\olww{\infty}e^{c_1|\psi_\eps-\sigma|}\leq C\ogkos{\qwrio}
\]
which implies that 
\[
\left(\olww{\infty} (u^2+\eps^2)^{\frac{c_1}{2}}\right)  \left(\olww{}
(u^2+\eps^2)^{-\frac{c_1}{2}}\right)\leq C^2 \ogkos{\qwrio}^2
\]
and for small \(p>0\):
\[
\mkf{\qwrio_\infty} u_\eps\geq
C\left(\frac{1}{v_\infty}\olww{\infty}u_\eps^p\right)^{1/p}
\] 
We follow again  \cite{lspde} appealing to \eqref{moiq}  and get the bound:
\[
\mkf{\qwrio_\infty} |u|\geq
C\left(\frac{1}{v_\infty}\olww{\infty}u^p\right)^{1/p}
\]
for any \(p\leq \frac{n}{n-2}\). 
\subsection{The two dimensional case}
We will derive a version of Hardy's inequality for the 
two dimensional situation that is not covered in the general case.
Therefore we start with the radial blow-up of the plane covered in the 
\[
C_1=\xw{2}\setminus 
\{ |x_1|\geq \eps |x_2|\} , C_2=\xw{2}\setminus \{|x_1|\geq
\eps |x_2|\} 
\]
We set in \(C_1\)
\[
x_1=\frac{r\xi}{\sqrt{1+\xi^2}},\quad x_2=\frac{r}{\sqrt{1+\xi^2}}
\] 
and interchange the roles of \(x_1,x_2\) in \(C_2\) We obtain for 
\[
P(\dian{x})= r^mR_j(\xi),\qquad j=1,2
\]
the elementary identity:
\[
|\kl P|^2=r^{2(m-1)}\left[\frac{m^2R_j^2}{(1+\xi^2)^m}
+(1+\xi^2)^2r^2\frac{d}{d\xi}
\left((1+\xi^2)^{-\frac{m}{2}} R_j\right)^2\right] 
\]
Therefore we have that:
\[
\frac{P^{2(1-\frac1m)}}{|\kl P|^2}=
\frac{R^2}{m^2+\left((1+\xi^2)R'-m\xi R\right)^2}\leq \frac{R^2}{m^2}
\]
Similarly:
\[
\left|\frac{\kl P}{P}\right|^2=\frac{m^2}{r^2}+
\left[(1+\xi^2)\frac{R'}{R}-m\xi\right]^2
\]
We compute then that for \(\delta<1\) and  
\(f\in C_0^\infty(\xw{2}\setminus \left(\{ P=0\}\bigcup \{\log
|P|=\delta\}\right)) \) then 
\[
\int_{\xw{2}} \frac{1}{P^{\frac2m}|\log\frac{|P|}{\delta}|}f^2\leq
\frac{C_1(P)}{|\log(\delta)|}\int_{\xw{2}} |\kl f|^2
\]
and
\[
\int_{\xw{2}} \frac{|\kl P|^2}{P^2|\log\frac{|P|}{\delta}|}f^2\leq 
\frac{C_2(P)}{|\log(\delta)|} \int_{\xw{2}} |\kl f|^2
\]
Then we start localizing in areas of constant sign for \(R\):
\[
I_{j,\eps}=
\int_{\xw{2}}\frac{1}{|\log\frac{|P|}{\delta}|}\left|\frac{\kl P}{P}\right|^2 
\chi_{j,\eps}^2f^2\leq 
C_\eps\int_{\xw{2}}\frac{m^2}{r|\log \frac{r}{\delta}|}\chi_jf^2+
\frac{r}{|\log\frac{r}{\delta}|}
\left((1+\xi^2)\frac{R'}{R}-m\xi\right)^2\chi_jf^2 
\]
We used the inequality:
\[
(a+b)^2\geq (1-\eps)^2\left(a^2-\frac{1}{\eps^2}b^2\right) 
\]
alternatively in the regions \((\log R)^2\geq \eta^2 (\log r)^2\) and 
\((\log R)^2\leq 2\eta^2(\log r)^2\).
 The inequality is majorised after integration by parts in the radial variable
through the elementary inequality:
\[
\int_0^\infty \frac{1}{r^2|\log\frac{r}{\delta}|} g^2 dr \leq
\frac{4}{|\log\delta|}\int_0^\infty g^2
\]
This proved easily by splitting the integral after arranging \(\eps\) 
according to the support of \(g\)
\[
\int_0^\infty \frac{1}{r^2|\log\frac{r}{\delta}|} g^2 dr=
\int_0^{\delta^{1+\eps}}  \frac{1}{r^2|\log\frac{r}{\delta}|} g^2 dr
+\int_{\delta^{-(1+\eps)}}^\infty  \frac{1}{r^2|\log\frac{r}{\delta}|} g^2 dr
\]
The first integral is written after integration by parts:
\[
2\int_0^{\delta^{1+\eps}}
\frac{1}{r\log {\delta}{r}}gg' +\frac{1}{r^2(\log\frac{\delta}{r})^2}g^2\leq 
\frac1\varepsilon\int_0^\infty g'2 + \frac{1+\varepsilon}{|\log\delta|}
\int_0^{\delta^{1+\eps}} \frac{1}{r^2|\log\frac{\delta}{r}|}f^2
\]
then we choose \(\varepsilon=\frac{|\log\delta|}{2}\) and we get:
\[
\int_0^{\delta^{1+\eps}}\frac{1}{r^2\log {\delta}{r}}g^2\leq 
\frac{C}{\eps|\log\delta|}\int_0^{\delta^{1+\eps}}g'2
\]
Similiarly for the other integral. 
setting as \(g^2=rf^2\) and splitting the integral in two pieces then
\[
\int_0^\infty  \frac{1}{r^2|\log \frac{r}{\delta}|} rf^2\leq 
\frac{2C_\eps}{|\log\delta|}
 \int_0^\infty  rf'^2+\frac{2C_\eps}{|\log\delta|}
\int_0^\infty \frac{1}{r^2|\log(\frac{r}{\delta})|} rf^2\leq
\]
\[
\leq \frac{2C_\eps}{|\log\delta|}\int_0^\infty f'^2rdr
\]
The two dimenional inequality is obtained by a direct application of the
usual one dimensional inequality in the  
\(\xi\)-variable after the formula:
\[
\left(\frac{R'}{R}\right)^2\leq 2\left(C+
\sum_{j=1}^m \frac{m_j^2}{(\xi-\xi_j)^2}\right)
\]
\subsection{Curvature estimates}
In the integration by parts formulas we substitute the curvature identities
 we find that:
\kseks
\begin{split}
\olws\phi^2|\kl\rie|^2\leq 3\olws \phi^2(|\R|^3+|\rie|^2),\qquad
\olws\phi^2|\kl\ric|^2\leq 3\olws \phi^2|\R|^3
\end{split}
\tleks

The iteration scheme suggests that
\[
\eaf{\qwrio^*}|\rie|\leq \hco(\eta,\mu)\mkf{\qwrio^*}|\R|
\]
\[
\eaf{\qwrio^*}|\ric|\leq \hco(\eta,\mu)\mkf{\qwrio^*}|\R|
\]

\subsection{Local properties of eigenfunctions}

In the case of an eigenfunction we have the following
\begin{subequations}
\begin{align}
\olws \phi^2|\kl^2 u_\ei|^2 & \leq \ei\left(\ei+
\hco(\eta,\mu)||\ric||_{1,\qwrio}\right)
\left(\olws \phi^2u_\ei^2\right)\label{eig1}\\
\olws \phi^2 |\kl^3u_\ei|^2 &\leq
\ei^2\left(\ei+\hco(\eta,\mu)||\ric||_{1,\qwrio}\right)
\olws \phi^2u_\ei^2\label{eig2}
\end{align}
\end{subequations}
Performing partial integration to the term:
\[
\olws \phi^2|\ric||\kl f|^2=-2\olws |\ric|f\phi\kl\phi\cdot \kl f-
\olws \phi^2 f\kl|\ric|\cdot\kl f-\olws |\ric|\phi^2f\lap{g} f
\]
Young's inequality along with harmonic approximation for 
\(\sqrt{|\ric|^2+\eps}\) leads to 
\[
\olws \phi^2 |\ric||\kl u_\ei|^2 \leq C\ei \olws |\ric|\phi^2u_\ei^2
\]
Similarly we get for 
Finally we conclude that
\begin{subequations}\label{eig}
\begin{align}
\olws \phi^2|\kl^2 u_\ei|^2 & \leq \ei^2\olws 
\phi^2\left(1+\frac{|\ric|}{\ei}\right)\label{eig1}\\
\olws \phi^2 |\kl^3u_\ei|^2 &\leq C\ei^3\olws \phi^2(1+\frac{|\rie|}{\ei})
u_\ei^2+\olws \left[\frac{1}{\eps^2}\left(|\kl\rie|^2+|\kl\ric|^2\right)+|\rie|^2+|\ric|^2\right]
\phi^2|\kl u_\ei|^2\label{eig2}
\end{align}
\end{subequations}

\subsubsection{Harnack inequalities}
\paragraph{The eigenfunction.}  We have for   \(\gamma=\ei, 
\tilde{u}_{\ei,\eps}=\sqrt{u_\ei^2+\eps^2} \) 
that
\kseks
\eaf{\qwrio} \tilde{u}_{\ei,\eps}\leq \hco(\eta,\ei)
\mkf{\qwrio}\tilde{u}_{\ei,\eps}
\tleks
\paragraph{The gradient. }
Now the gradient  \(G_{\ei,\eps}=|\kl u_\ei|^2 +\eps^2\) 
requires that we use the  \(\gamma_1=\ei+\kappa \) and we get that
 \kseks
\eaf{\qwrio} G_{\ei,\eps} \leq
\hco(\eta,\ei)\left(\ei+\hco(\eta,\mu)||\ric||_{1,\qwrio}\right)
\mkf{\qwrio}G_{\ei,\eps}
\tleks
\paragraph{The hessian estimate}
The estimate
\begin{eqnarray}
\eaf{\qwrio}H_{\ei,\eps} 
\leq \hco(\eta,\ei)\ei^2\left(\ei+\hco(\eta,\mu)||\ric||_{1,\qwrio}\right)
\mkf{\qwrio}H_{\ei,\eps}
\end{eqnarray}
\paragraph{The estimates for the restriction on the spherical front} The 
restriction of the eigenfunction 
\[
u(r,\dian{\theta})=e^{\duna}\sin(\beta_1\phi)
\]
of the spherical front satifies the following inequalities for \(j=1,2\):
\[
\olwf{} \theta^2| \klgs^ju|^2\leq C\olwf{}\theta^2(\ei+R^2)u^2
\]
where \(R\) is a polynomial depending on 
\(\rie,\kl\rie,\dots,\kl^2\rie,\ric,\kl\ric,\kl\ric\). This combined with
Michael-Simon Sobolev inequality provides Harnack inequalities for the 
restriction of \(u,|\klgs u|,\klgs^2u\) on the spherical front. 
\subsubsection{The Bernstein inequalities}

The integration by parts formulas suggest along with the Harnack
inequalities the following Berstein type
inequalities in geodesic pixels:

\begin{berniq}

The following estimates hold in a domain inside a 
geodesic pixel \(\qwrio\subset \bricks\)

\begin{align}
|\kl^2u_\ei| & \leq \eaf{\qwrio}|\kl^2u_\ei|\leq
C_2(\qwrio) \ei^{\frac{n}{2}+1}\left(|\kl u_\ei|+\eps\right) \leq 
C_3(\qwrio) \ei^{\frac{n}{2}+2}\left(|u_\ei|+\eps\right)  \\
|\kl u_\ei| &\leq \eaf{\qwrio}|\kl u_\ei|\leq
C_4(\qwrio)\ei^{\frac{n}{2}+1}(|u_\ei|+\eps)
\end{align}

\end{berniq}

\begin{flushleft}
\noindent 
Center for Plasma Physics and Laser, \\ 
Department of Electronics,\\
TEI Crete,\,Gr-73103 \\
email: {\ttfamily dpliakis@gmail.com}
\end{flushleft}
\end{document}